\documentclass{article}
\usepackage{amsmath}

\input{tcilatex}
\begin{document}

\date{}
\author{Demetrios Serakos$^{1}$, John E. Gray$^{2}$ and Hazim Youssef$^{1}$ 
\\
Naval Surface Warfare Center\\
$^{1}$Warfare Systems Department\\
$^{2}$Electromagnetic and Sensor Systems Department\\
Dahlgren, VA 22448}
\title{TOPICS\ IN\ MITIGATING RADAR\ BIAS}
\maketitle

\begin{abstract}
In this paper, we investigate two topics related to mitigating the effect of
radar bias in ballistic missile tracking applications. We determine the
absolute bias between two radars in polar coordinates when their relative
bias is given in rectangular coordinates. Using this result, we then obtain
the optimized steady-state filter to handle the random bias.
\end{abstract}

\section{\protect\bigskip Introduction}

There are several facets to the problem of tracking ballistic missiles with
radar that require enhanced error correction to effectively track threats.
In this paper, we obtain the exact form of the bias error for the coordinate
transformation problem. This result is useful in Ballistic Missile Defense
bistatic applications where one sensor is used for launching an interceptor,
while another is used to track the threat. Thus, the problem of translation
between internal sensor coordinate frames to a common frame (that is, used
by all sensors) is important. The coordinate transformation problem from
Cartesian to spherical coordinates introduces a bias that, if accounted for,
can be corrected in the design of a filter. This problem occurs when one has
multiple launch platforms, because each local track must be formatted for a
common reference frame. When bias correction is accomplished correctly, one
can improve tracking performance of the filter and increase the likelihood
that an interceptor can successfully engage a threat.

\section{An Optimized Method of Obtaining Absolute Bias}

Although relative bias calculation can be used to provide correct
association of tracks from two sensors, the calculation of the absolute bias
is required to correct the track state and is needed for track fusion and
for producing a Single Integrated Air Picture. Methods for obtaining the
relative bias between two radars tracking the same ballistic missile are
presented in Levedahl \cite{Mark} and Brown, Weisman and Brock \cite{BWB}.
The methods presented in these reports have to do with maximizing a
likelihood function. The relative biases obtained in these papers are
determined in rectangular coordinates. In this paper, the absolute bias for
the two sensors is calculated from the relative bias by solving a
minimization problem. The problem is set up to minimize the weighted sum of
the two absolute biases while viewing the given relative bias as a
constraint.

A point in $3$-dimensional space in both rectangular and spherical
coordinates\footnote{%
Denote yaw (azimuth) by $\psi $, pitch (elevation) by $\theta $. $\phi $ is
normally reserved for roll; however, roll is not used here.} is denoted by:%
\begin{equation}
\overrightarrow{p}=\left[ 
\begin{array}{c}
x \\ 
y \\ 
z%
\end{array}%
\right] ;\text{ and}\qquad \overrightarrow{\pi }=\left[ 
\begin{array}{c}
r \\ 
\psi \\ 
\theta%
\end{array}%
\right] \text{, respectively.}  \label{1a}
\end{equation}

The transformations between the coordinates are $\overrightarrow{p}=f(%
\overrightarrow{\pi })$ and $\overrightarrow{\pi }=f^{-1}(\overrightarrow{p}%
) $, which are given by:%
\begin{equation}
f(\overrightarrow{\pi })=\left[ 
\begin{array}{c}
r\cos \theta \cos \psi \\ 
r\cos \theta \sin \psi \\ 
r\sin \theta%
\end{array}%
\right] ;\qquad \qquad f^{-1}(\overrightarrow{p})=\left[ 
\begin{array}{c}
\sqrt{x^{2}+y^{2}+z^{2}} \\ 
\arctan \left( y/x\right) \\ 
\arctan \left( z/\sqrt{x^{2}+y^{2}}\right)%
\end{array}%
\right] \text{ .}  \label{1b}
\end{equation}%
We need the following definitions:

\begin{center}
\begin{tabular}{||l|l||}
\hline\hline
$P_{1}$ & Target position as seen by sensor 1 \\ \hline
$P_{2}$ & Target position as seen by sensor 2 \\ \hline
$P_{T}$ & True target position (unknown) \\ \hline
$B_{1}$ & Sensor 1 bias \\ \hline
$B_{2}$ & Sensor 2 bias \\ \hline
$B_{R}$ & Relative bias \\ \hline
$P_{1TO2}$ & Sensor 2 position from sensor 1 \\ \hline\hline
\end{tabular}
\end{center}

$P_{1,ENU(1)}=(x_{1},y_{1},z_{1})_{ENU(1)}^{\prime }$, $%
P_{2,ENU(2)}=(x_{2},y_{2},z_{2})_{ENU(2)}^{\prime }$. $B_{1,ENU(1)}=(\Delta
x_{1},\Delta y_{1},\Delta z_{1})_{ENU(1)}^{\prime }$, $B_{2,ENU(2)}=(\Delta
x_{2},\Delta y_{2},\Delta z_{2})_{ENU(2)}^{\prime }$. (ENU denotes the East
North Up coordinate system.) Thus, we have in the sensor coordinates%
\begin{equation}
P_{T,ENU(1)}=P_{1,ENU(1)}+B_{1,ENU(1)}  \label{1c}
\end{equation}%
\begin{equation}
P_{T,ENU(2)}=P_{2,ENU(2)}+B_{2,ENU(2)}\text{ .}  \label{1d}
\end{equation}%
If we use an ENU coordinate system located at sensor 1, (\ref{1d}) becomes%
\begin{equation}
P_{T,ENU(1)}=P_{1TO2,ENU(1)}+P_{2,ENU(1)}+B_{2,ENU(1)}  \label{1e}
\end{equation}%
where $P_{1TO2,ENU(1)}$ is the position vector from the first sensor to the
second sensor in $ENU(1)$. The relative bias in $ENU(1)$ is 
\begin{equation*}
B_{R,ENU(1)}=B_{2,ENU(1)}-B_{1,ENU(1)}
\end{equation*}%
\begin{equation*}
=(P_{T,ENU(1)}-P_{1TO2,ENU(1)}-P_{2,ENU(1)})-(P_{T,ENU(1)}-P_{1,ENU(1)})
\end{equation*}%
\begin{equation}
=P_{1,ENU(1)}-P_{1TO2,ENU(1)}-P_{2,ENU(1)}\text{ .}  \label{1f}
\end{equation}%
%
%
%
%
%
%
%
%
%
%
%
%
%
%
%
%
%
%
%
%
%
%
%
%
%
%
%
%
%
%
%
%
%
%
%
%
%
%
%
%
%
%
%

We consider the coordinate transformations to allow us to go from $ENU$ to
radar-face coordinates for a particular sensor. Each sensor has its own face
and ENU coordinate systems. The face coordinate system (denoted FACE) of a
sensor is related to the ENU coordinate system of a sensor by the following
transformation:%
\begin{equation*}
T_{ENU(i)2FACE(i)}=\ \left[ 
\begin{array}{ccc}
\cos \theta _{i} & 0 & \sin \theta _{i} \\ 
0 & 1 & 0 \\ 
-\sin \theta _{i} & 0 & \cos \theta _{i}%
\end{array}%
\right] \left[ 
\begin{array}{ccc}
\cos \psi _{i} & \sin \psi _{i} & 0 \\ 
-\sin \psi _{i} & \cos \psi _{i} & 0 \\ 
0 & 0 & 1%
\end{array}%
\right]
\end{equation*}

\begin{equation}
=\left[ 
\begin{array}{ccc}
\cos \theta _{i}\cos \psi _{i} & \cos \theta _{i}\sin \psi _{i} & \sin
\theta _{i} \\ 
-\sin \psi _{i} & \cos \psi _{i} & 0 \\ 
-\sin \theta _{i}\cos \psi _{i} & -\sin \theta _{i}\sin \psi _{i} & \cos
\theta _{i}%
\end{array}%
\right] \allowbreak  \label{1i}
\end{equation}%
where $i=1,2$. We also have that%
\begin{equation}
T_{FACE(i)2ENU(i)}=\left[ 
\begin{array}{ccc}
\cos \theta _{i}\cos \psi _{i} & -\sin \psi _{i} & -\sin \theta _{i}\cos
\psi _{i} \\ 
\cos \theta _{i}\sin \psi _{i} & \cos \psi _{i} & -\sin \theta _{i}\sin \psi
_{i} \\ 
\sin \theta _{i} & 0 & \cos \theta _{i}%
\end{array}%
\right] \text{ ,}\allowbreak  \label{1j}
\end{equation}%
which is the transpose of (\ref{1i}). We can also have the matrix $%
T_{ENU(i)2FACE(j)}$, which is%
\begin{equation}
T_{ENU(i)2FACE(j)}=\left[ 
\begin{array}{ccc}
\cos \theta _{i,j}\cos \psi _{i,j} & \cos \theta _{i,j}\sin \psi _{i,j} & 
\sin \theta _{i,j} \\ 
-\sin \psi _{i,j} & \cos \psi _{i,j} & 0 \\ 
-\sin \theta _{i,j}\cos \psi _{i,j} & -\sin \theta _{i,j}\sin \psi _{i,j} & 
\cos \theta _{i,j}%
\end{array}%
\right] \text{ .}  \label{1k}
\end{equation}%
The absolute (as opposed to relative) bias can be expressed in the face
coordinates:%
\begin{equation}
B_{i,FACE(i)}=\Delta r\cdot \overrightarrow{u_{r}}+\Delta c_{A}\cdot 
\overrightarrow{u_{cA}}+\Delta c_{B}\cdot \overrightarrow{u_{cB}}  \label{1l}
\end{equation}%
where $\overrightarrow{u_{r}}$ is the unit vector in the range coordinate
and $\overrightarrow{u_{cA}},$ $\overrightarrow{u_{cB}}$ are the two cross
range coordinate unit vectors.\ Substituting $p_{T}\Delta \psi =\Delta c_{A}$
and $p_{T}\Delta \theta =\Delta c_{B}$ where $p_{Ti}=\left\Vert
P_{T}(i)\right\Vert $ (see note\footnote{%
True position is not available. When applying this method measured position
is used for this calculation instead.}), the distance from sensor to the
target, we get%
\begin{equation}
B_{i,FACE(i)}=\Delta r_{i}\cdot \overrightarrow{u_{r}}+p_{Ti}\Delta \psi
_{i}\cdot \overrightarrow{u_{cA}}+p_{Ti}\Delta \theta _{i}\cdot 
\overrightarrow{u_{cB}}  \label{1m}
\end{equation}%
\begin{equation*}
B_{i,ENU(i)}=\left[ 
\begin{array}{ccc}
\cos \theta _{i}\cos \psi _{i} & -\sin \psi _{i} & -\sin \theta _{i}\cos
\psi _{i} \\ 
\cos \theta _{i}\sin \psi _{i} & \cos \psi _{i} & -\sin \theta _{i}\sin \psi
_{i} \\ 
\sin \theta _{i} & 0 & \cos \theta _{i}%
\end{array}%
\right] \left[ 
\begin{array}{c}
\Delta r_{i} \\ 
p_{Ti}\Delta \psi _{i} \\ 
p_{Ti}\Delta \theta _{i}%
\end{array}%
\right]
\end{equation*}%
\begin{equation}
=\left[ 
\begin{array}{c}
\Delta r_{i}\cos \theta _{i}\cos \psi _{i}-\Delta \psi _{i}\left( \sin \psi
_{i}\right) \cdot p_{Ti}-\Delta \theta _{i}\left( \sin \theta _{i}\cos \psi
_{i}\right) \cdot p_{Ti} \\ 
\Delta r_{i}\cos \theta _{i}\sin \psi _{i}+\Delta \psi _{i}\left( \cos \psi
_{i}\right) \cdot p_{Ti}-\Delta \theta _{i}\left( \sin \theta _{i}\sin \psi
_{i}\right) \cdot p_{Ti} \\ 
\Delta r_{i}\sin \theta _{i}+\Delta \theta _{i}\left( \cos \theta
_{i}\right) \cdot p_{Ti}%
\end{array}%
\right] \text{ .}  \label{1n}
\end{equation}%
We can obviously obtain $B_{i,ENU(j)}$ (for $i$ not necessarily equal to $j$%
) if needed. The quantities $\Delta r_{1},$ $\Delta \psi _{1},$ $\Delta
\theta _{1},$ $\Delta r_{2},$ $\Delta \psi _{2},$ $\Delta \theta _{2}$ are
the ones we minimize. To refer to these as a group we on occasion write $%
e=\left( \Delta r_{1},\Delta \psi _{1},\Delta \theta _{1},\Delta
r_{2},\Delta \psi _{2},\Delta \theta _{2}\right) $.

We need tolerances or costs for the sensor biases. These are expressed in
spherical coordinates.

\begin{center}
\begin{tabular}{||l|l||}
\hline\hline
$k_{r1}$ & Sensor 1 range bias cost, unitless \\ \hline
$k_{r2}$ & Sensor 2 range bias cost, unitless \\ \hline
$k_{\psi 1}$ & Sensor 1 azimuth bias cost, meters \\ \hline
$k_{\psi 2}$ & Sensor 2 azimuth bias cost, meters \\ \hline
$k_{\theta 1}$ & Sensor 1 elevation bias cost, meters \\ \hline
$k_{\theta 2}$ & Sensor 2 elevation bias cost, meters \\ \hline\hline
\end{tabular}
\end{center}

\subsection{Problem Statement}

We want to compute the minimum (absolute) bias cost for the two sensors when
there are known (computed) expressions for the relative bias. The given
relative bias is expressed in ENU rectangular coordinates. We compute the
minimum absolute bias in spherical coordinates. The relative bias in
rectangular coordinates contrasted with the absolute bias in spherical
coordinates allows us to formulate this as a minimization problem. We view
the relative bias as a constraint. We use a quadratic cost: 
\begin{equation}
F=\frac{k_{r_{1}}^{2}}{2}\cdot \left( \Delta r_{1}\right) ^{2}+\frac{k_{\psi
_{1}}^{2}}{2}\cdot \left( \Delta \phi _{1}\right) ^{2}+\frac{k_{\theta
_{1}}^{2}}{2}\cdot \left( \Delta \theta _{1}\right) ^{2}+\frac{k_{r_{2}}^{2}%
}{2}\cdot \left( \Delta r_{2}\right) ^{2}+\frac{k_{\psi _{2}}^{2}}{2}\cdot
\left( \Delta \phi _{2}\right) ^{2}+\frac{k_{\theta _{2}}^{2}}{2}\cdot
\left( \Delta \theta _{2}\right) ^{2}\text{ .}  \label{1o}
\end{equation}%
So that the addition in (\ref{1o}) is permissible, we have that $k_{r_{1}}$, 
$k_{r_{2}}$ are unitless and $k_{\psi _{1}}$, $k_{\theta _{1}}$, $k_{\psi
_{2}}$, $k_{\theta _{2}}$ are in meters. We note $F$ may be rewritten in the
form%
\begin{equation*}
F=\left[ 
\begin{array}{ccc}
\Delta r_{1} & \Delta \psi _{1} & \Delta \theta _{1}%
\end{array}%
\right] \left[ 
\begin{array}{ccc}
2/k_{r_{1}}^{2} & 0 & 0 \\ 
0 & 2/k_{\psi _{1}}^{2} & 0 \\ 
0 & 0 & 2/k_{\theta _{1}}^{2}%
\end{array}%
\right] ^{-1}\left[ 
\begin{array}{c}
\Delta r_{1} \\ 
\Delta \psi _{1} \\ 
\Delta \theta _{1}%
\end{array}%
\right]
\end{equation*}%
\begin{equation}
+\left[ 
\begin{array}{ccc}
\Delta r_{2} & \Delta \psi _{2} & \Delta \theta _{2}%
\end{array}%
\right] \left[ 
\begin{array}{ccc}
2/k_{r_{2}}^{2} & 0 & 0 \\ 
0 & 2/k_{\psi _{2}}^{2} & 0 \\ 
0 & 0 & 2/k_{\theta _{2}}^{2}%
\end{array}%
\right] ^{-1}\left[ 
\begin{array}{c}
\Delta r_{2} \\ 
\Delta \psi _{2} \\ 
\Delta \theta _{2}%
\end{array}%
\right] \text{ ,}  \label{1p}
\end{equation}%
which we recognize as being in the form of a Mahalanobis distance. We note
that the Mahalanobis distance comes up in the Levedahl/Lincoln Labs work (%
\cite{Mark} and \cite{BWB}) when the log is taken of the Gaussian
distribution. The cost $F$ is minimized subject to this equality constraint:%
\begin{equation}
G(B)=\left( B_{2,ENU(1)}-B_{1,ENU(1)}\right) -B_{R,ENU(1)}=0\text{ .}
\label{1p1}
\end{equation}%
Thus, we have%
\begin{equation*}
G(B)=\left[ 
\begin{array}{c}
\Delta r_{2}\cos \theta _{2}\cos \psi _{2}-\Delta \psi _{2}\left( \sin \psi
_{2}\right) \cdot p_{T2}-\Delta \theta _{2}\left( \sin \theta _{2}\cos \psi
_{2}\right) \cdot p_{T2} \\ 
\Delta r_{2}\cos \theta _{2}\sin \psi _{2}+\Delta \psi _{2}\left( \cos \psi
_{2}\right) \cdot p_{T2}-\Delta \theta _{2}\left( \sin \theta _{2}\sin \psi
_{2}\right) \cdot p_{T2} \\ 
\Delta r_{2}\sin \theta _{2}+\Delta \theta _{2}\left( \cos \theta
_{2}\right) \cdot p_{T2}%
\end{array}%
\right]
\end{equation*}%
\begin{equation}
-\left[ 
\begin{array}{c}
\Delta r_{1}\cos \theta _{1}\cos \psi _{1}-\Delta \psi _{1}\left( \sin \psi
_{1}\right) \cdot p_{T1}-\Delta \theta _{1}\left( \sin \theta _{1}\cos \psi
_{1}\right) \cdot p_{T1} \\ 
\Delta r_{1}\cos \theta _{1}\sin \psi _{1}+\Delta \psi _{1}\left( \cos \psi
_{1}\right) \cdot p_{T1}-\Delta \theta _{1}\left( \sin \theta _{1}\sin \psi
_{1}\right) \cdot p_{T1} \\ 
\Delta r_{1}\sin \theta _{1}+\Delta \theta _{1}\left( \cos \theta
_{1}\right) \cdot p_{T1}%
\end{array}%
\right] -B_{R}=0  \label{1q}
\end{equation}%
where all of the terms in (\ref{1q}) reside entirely in one or the other of
the two $ENU$ coordinate systems. We see that $G(B)$ gives that the
difference between the two absolute biases (whatever they may be) is equal
to the relative bias. Also, we note that (\ref{1q}) is affine. Another
equivalent representation for $G(B)$ is%
\begin{equation}
G(B)=A\left( p_{T2},\psi _{2},\theta _{2}\right) \left[ 
\begin{array}{c}
\Delta r_{2} \\ 
\Delta \psi _{2} \\ 
\Delta \theta _{2}%
\end{array}%
\right] -A\left( p_{T1},\psi _{1},\theta _{1}\right) \left[ 
\begin{array}{c}
\Delta r_{1} \\ 
\Delta \psi _{1} \\ 
\Delta \theta _{1}%
\end{array}%
\right] -B_{R}  \label{1s}
\end{equation}%
where%
\begin{equation}
A\left( p_{T1},\psi _{1},\theta _{1}\right) =\left[ 
\begin{array}{ccc}
\cos \theta _{1}\cos \psi _{1} & -\sin \psi _{1}\cdot p_{T1} & -\sin \theta
_{1}\cos \psi _{1}\cdot p_{T1} \\ 
\cos \theta _{1}\sin \psi _{1} & \cos \psi _{1}\cdot p_{T1} & -\sin \theta
_{1}\sin \psi _{1}\cdot p_{T1} \\ 
\sin \theta _{1} & 0 & \cos \theta _{1}\cdot p_{T1}%
\end{array}%
\right]  \label{1t}
\end{equation}%
,%
\begin{equation}
A\left( p_{T2},\psi _{2},\theta _{2}\right) =\left[ 
\begin{array}{ccc}
\cos \theta _{2}\cos \psi _{2} & -\sin \psi _{2}\cdot p_{T2} & -\sin \theta
_{2}\cos \psi _{2}\cdot p_{T2} \\ 
\cos \theta _{2}\sin \psi _{2} & \cos \psi _{2}\cdot p_{T2} & -\sin \theta
_{2}\sin \psi _{2}\cdot p_{T2} \\ 
\sin \theta _{2} & 0 & \cos \theta _{2}\cdot p_{T2}%
\end{array}%
\right] \text{ .}  \label{1u}
\end{equation}%
Setting $G(B)=0$, we solve for $\Delta r_{2},\Delta \psi _{2},\Delta \theta
_{2}$%
\begin{equation}
\left[ 
\begin{array}{c}
\Delta r_{2} \\ 
\Delta \psi _{2} \\ 
\Delta \theta _{2}%
\end{array}%
\right] =A^{-1}\left( p_{T2},\psi _{2},\theta _{2}\right) \left( A\left(
p_{T1},\psi _{1},\theta _{1}\right) \left[ 
\begin{array}{c}
\Delta r_{1} \\ 
\Delta \psi _{1} \\ 
\Delta \theta _{1}%
\end{array}%
\right] +B_{R}\right)  \label{1v}
\end{equation}%
%
%
%
%
%
%
%
%
%
%
%
%
%
%
%
%
%
%
%
%
%
%
%
%
%
%
%
%
%
%
%
%
%
%
%
%
%
%
%
%
%
%
%
(provided $p_{T2}\neq 0$.) This vector equality constraint (\ref{1q}) can be
written in the form of three scalar equality constraints%
\begin{equation*}
G_{E}(B)=\Delta r_{2}\cos \theta _{2}\cos \psi _{2}-\Delta \psi _{2}\sin
\psi _{2}\cdot p_{T2}-\Delta \theta _{2}\sin \theta _{2}\cos \psi _{2}\cdot
p_{T2}
\end{equation*}%
\begin{equation}
-\Delta r_{1}\cos \theta _{1}\cos \psi _{1}+\Delta \psi _{1}\sin \psi
_{1}\cdot p_{T1}+\Delta \theta _{1}\sin \theta _{1}\cos \psi _{1}\cdot
p_{T1}-B_{RE}  \label{1w}
\end{equation}%
\begin{equation*}
G_{N}(B)=\Delta r_{2}\cos \theta _{2}\sin \psi _{2}+\Delta \psi _{2}\cos
\psi _{2}\cdot p_{T2}-\Delta \theta _{2}\sin \theta _{2}\sin \psi _{2}\cdot
p_{T2}
\end{equation*}%
\begin{equation}
-\Delta r_{1}\cos \theta _{1}\sin \psi _{1}-\Delta \psi _{1}\cos \psi
_{1}\cdot p_{T1}+\Delta \theta _{1}\sin \theta _{1}\sin \psi _{1}\cdot
p_{T1}-B_{RN}  \label{1x}
\end{equation}%
\begin{equation}
G_{U}(B)=\Delta r_{2}\sin \theta _{2}+\Delta \theta _{2}\cos \theta
_{2}\cdot p_{T2}-\Delta r_{1}\sin \theta _{1}-\Delta \theta _{1}\cos \theta
_{1}\cdot p_{T1}-B_{RU}\text{ .}  \label{1y}
\end{equation}%
%
%
%
%
%
%
%
%
%
%
%
%
%
%
%
%
%
%
%
%
%
%
%
%
%
%
%
%
%
%
%
%
%
%
%
%
%
%
%
%
%
%
%

\subsection{Solving The Minimization Problem}

To solve this minimization problem, we need to take a few derivatives. We
need the gradient of the function to be minimized. We also need the gradient
of the constraint, which is an equality constraint in this case.%
\begin{equation}
\nabla F=\left[ 
\begin{array}{c}
\partial F/\partial \Delta r_{1} \\ 
\partial F/\partial \Delta \psi _{1} \\ 
\partial F/\partial \Delta \theta _{1} \\ 
\partial F/\partial \Delta r_{2} \\ 
\partial F/\partial \Delta \psi _{2} \\ 
\partial F/\partial \Delta \theta _{2}%
\end{array}%
\right] =\left[ 
\begin{array}{c}
k_{r_{1}}^{2}\cdot \Delta r_{1} \\ 
k_{\psi _{1}}^{2}\cdot \Delta \psi _{1} \\ 
k_{\theta _{1}}^{2}\cdot \Delta \theta _{1} \\ 
k_{r_{2}}^{2}\cdot \Delta r_{2} \\ 
k_{\psi _{2}}^{2}\cdot \Delta \psi _{2} \\ 
k_{\theta _{2}}^{2}\cdot \Delta \theta _{2}%
\end{array}%
\right]   \label{1z}
\end{equation}%
\begin{equation}
\nabla G_{E}=\left[ 
\begin{array}{c}
\partial G_{E}/\partial \Delta r_{1} \\ 
\partial G_{E}/\partial \Delta \psi _{1} \\ 
\partial G_{E}/\partial \Delta \theta _{1} \\ 
\partial G_{E}/\partial \Delta r_{2} \\ 
\partial G_{E}/\partial \Delta \psi _{2} \\ 
\partial G_{E}/\partial \Delta \theta _{2}%
\end{array}%
\right] =\left[ 
\begin{array}{c}
-\cos \theta _{1}\cos \psi _{1} \\ 
\sin \psi _{1}\cdot p_{T1} \\ 
\sin \theta _{1}\cos \psi _{1}\cdot p_{T1} \\ 
\cos \theta _{2}\cos \psi _{2} \\ 
-\sin \psi _{2}\cdot p_{T2} \\ 
-\sin \theta _{2}\cos \psi _{2}\cdot p_{T2}%
\end{array}%
\right]   \label{1aa}
\end{equation}%
\begin{equation}
\nabla G_{N}=\left[ 
\begin{array}{c}
\partial G_{N}/\partial \Delta r_{1} \\ 
\partial G_{N}/\partial \Delta \psi _{1} \\ 
\partial G_{N}/\partial \Delta \theta _{1} \\ 
\partial G_{N}/\partial \Delta r_{2} \\ 
\partial G_{N}/\partial \Delta \psi _{2} \\ 
\partial G_{N}/\partial \Delta \theta _{2}%
\end{array}%
\right] =\left[ 
\begin{array}{c}
-\cos \theta _{1}\sin \psi _{1} \\ 
-\cos \psi _{1}\cdot p_{T1} \\ 
\sin \theta _{1}\sin \psi _{1}\cdot p_{T1} \\ 
\cos \theta _{2}\sin \psi _{2} \\ 
\cos \psi _{2}\cdot p_{T2} \\ 
-\sin \theta _{2}\sin \psi _{2}\cdot p_{T2}%
\end{array}%
\right]   \label{1ab}
\end{equation}%
\begin{equation}
\nabla G_{U}=\left[ 
\begin{array}{c}
\partial G_{U}/\partial \Delta r_{1} \\ 
\partial G_{U}/\partial \Delta \psi _{1} \\ 
\partial G_{U}/\partial \Delta \theta _{1} \\ 
\partial G_{U}/\partial \Delta r_{2} \\ 
\partial G_{U}/\partial \Delta \psi _{2} \\ 
\partial G_{u}/\partial \Delta \theta _{2}%
\end{array}%
\right] =\left[ 
\begin{array}{c}
-\sin \theta _{1} \\ 
0 \\ 
-p_{T1}\cdot \cos \theta _{1} \\ 
\sin \theta _{2} \\ 
0 \\ 
p_{T2}\cdot \cos \theta _{2}%
\end{array}%
\right] \text{ .}  \label{1ac}
\end{equation}%
We are looking for an optimal solution located at the point $e^{\ast
}=\left( \Delta r_{1}^{\ast },\Delta \psi _{1}^{\ast },\Delta \theta
_{1}^{\ast },\Delta r_{2}^{\ast },\Delta \psi _{2}^{\ast },\Delta \theta
_{2}^{\ast }\right) $. We employ the Kuhn-Tucker conditions that stipulate
the optimal solution $e^{\ast }$ should satisfy these equality constraints
for $e$ and there exist numbers $a_{1}^{\ast },a_{2}^{\ast },a_{3}^{\ast }$
such that%
\begin{equation}
\nabla F\left( e^{\ast }\right) =a_{1}^{\ast }\cdot \nabla G_{E}\left(
e^{\ast }\right) +a_{2}^{\ast }\cdot \nabla G_{N}\left( e^{\ast }\right)
+a_{3}^{\ast }\cdot \nabla G_{U}\left( e^{\ast }\right) \text{ .}
\label{1ad}
\end{equation}%
The gradients $\nabla G_{E}\left( e^{\ast }\right) ,\nabla G_{N}\left(
e^{\ast }\right) ,\nabla G_{U}\left( e^{\ast }\right) $ are linearly
independent. Taking an inventory of the equations and unknowns, we see that
there are $9$ unknowns ($e,a_{1},a_{2},a_{3}$) and $9$ equations ($3$ from
the equality constraint and $6$ from the above equation). We \textit{may} be
able to find the solution. Since the cost $F$ is quadratic and the
constraint $G$ is affine, the necessary conditions we give for optimality
are also sufficient conditions and an optimal solution $e^{\ast }$ is a
global optimal solution. Equation (\ref{1ad}) in longhand is:%
\begin{equation}
\left[ 
\begin{array}{c}
k_{r_{1}}^{2}\cdot \Delta r_{1} \\ 
k_{\psi _{1}}^{2}\cdot \Delta \psi _{1} \\ 
k_{\theta _{1}}^{2}\cdot \Delta \theta _{1} \\ 
k_{r_{2}}^{2}\cdot \Delta r_{2} \\ 
k_{\psi _{2}}^{2}\cdot \Delta \psi _{2} \\ 
k_{\theta _{2}}^{2}\cdot \Delta \theta _{2}%
\end{array}%
\right] =a_{1}\left[ 
\begin{array}{c}
-\cos \theta _{1}\cos \psi _{1} \\ 
\sin \psi _{1}\cdot p_{T1} \\ 
\sin \theta _{1}\cos \psi _{1}\cdot p_{T1} \\ 
\cos \theta _{2}\cos \psi _{2} \\ 
-\sin \psi _{2}\cdot p_{T2} \\ 
-\sin \theta _{2}\cos \psi _{2}\cdot p_{T2}%
\end{array}%
\right] +a_{2}\left[ 
\begin{array}{c}
-\cos \theta _{1}\sin \psi _{1} \\ 
-\cos \psi _{1}\cdot p_{T1} \\ 
\sin \theta _{1}\sin \psi _{1}\cdot p_{T1} \\ 
\cos \theta _{2}\sin \psi _{2} \\ 
\cos \psi _{2}\cdot p_{T2} \\ 
-\sin \theta _{2}\sin \psi _{2}\cdot p_{T2}%
\end{array}%
\right] +a_{3}\left[ 
\begin{array}{c}
-\sin \theta _{1} \\ 
0 \\ 
-p_{T1}\cdot \cos \theta _{1} \\ 
\sin \theta _{2} \\ 
0 \\ 
p_{T2}\cdot \cos \theta _{2}%
\end{array}%
\right] \text{ .}  \label{1ad1}
\end{equation}%
The right hand side of (\ref{1ad1}) may be written in the form of the
product of two matrices $M_{1}$ and $M_{2}$,%
\begin{equation}
\left[ 
\begin{array}{ccc}
-\cos \theta _{1}\cos \psi _{1} & -\cos \theta _{1}\sin \psi _{1} & -\sin
\theta _{1} \\ 
\sin \psi _{1}\cdot p_{T1} & -\cos \psi _{1}\cdot p_{T1} & 0 \\ 
\sin \theta _{1}\cos \psi _{1}\cdot p_{T1} & \sin \theta _{1}\sin \psi
_{1}\cdot p_{T1} & -p_{T1}\cdot \cos \theta _{1} \\ 
\cos \theta _{2}\cos \psi _{2} & \cos \theta _{2}\sin \psi _{2} & \sin
\theta _{2} \\ 
-\sin \psi _{2}\cdot p_{T2} & \cos \psi _{2}\cdot p_{T2} & 0 \\ 
-\sin \theta _{2}\cos \psi _{2}\cdot p_{T2} & -\sin \theta _{2}\sin \psi
_{2}\cdot p_{T2} & p_{T2}\cdot \cos \theta _{2}%
\end{array}%
\right] \left[ 
\begin{array}{c}
a_{1} \\ 
a_{2} \\ 
a_{3}%
\end{array}%
\right] =\left[ \QDATOP{M_{1}}{M_{2}}\right] \left[ 
\begin{array}{c}
a_{1} \\ 
a_{2} \\ 
a_{3}%
\end{array}%
\right] \text{ ,}  \label{1ae}
\end{equation}%
where%
\begin{equation}
M_{1}=\left[ 
\begin{array}{ccc}
-\cos \theta _{1}\cos \psi _{1} & -\cos \theta _{1}\sin \psi _{1} & -\sin
\theta _{1} \\ 
\sin \psi _{1}\cdot p_{T1} & -\cos \psi _{1}\cdot p_{T1} & 0 \\ 
\sin \theta _{1}\cos \psi _{1}\cdot p_{T1} & \sin \theta _{1}\sin \psi
_{1}\cdot p_{T1} & -p_{T1}\cdot \cos \theta _{1}%
\end{array}%
\right]   \label{1af}
\end{equation}%
\begin{equation}
M_{2}=\left[ 
\begin{array}{ccc}
\cos \theta _{2}\cos \psi _{2} & \cos \theta _{2}\sin \psi _{2} & \sin
\theta _{2} \\ 
-\sin \psi _{2}\cdot p_{T2} & \cos \psi _{2}\cdot p_{T2} & 0 \\ 
-\sin \theta _{2}\cos \psi _{2}\cdot p_{T2} & -\sin \theta _{2}\sin \psi
_{2}\cdot p_{T2} & p_{T2}\cdot \cos \theta _{2}%
\end{array}%
\right] \text{ .}  \label{1ag}
\end{equation}%
Note that $a_{1}$, $a_{2}$, $a_{3}$ have the units of meters. Let%
\begin{equation}
D_{1}=\left[ 
\begin{array}{ccc}
k_{r_{1}}^{2} & 0 & 0 \\ 
0 & k_{\psi _{1}}^{2} & 0 \\ 
0 & 0 & k_{\theta _{1}}^{2}%
\end{array}%
\right]   \label{1ah}
\end{equation}%
\begin{equation}
D_{2}=\left[ 
\begin{array}{ccc}
k_{r_{2}}^{2} & 0 & 0 \\ 
0 & k_{\psi _{2}}^{2} & 0 \\ 
0 & 0 & k_{\theta _{2}}^{2}%
\end{array}%
\right] \text{ .}  \label{1ai}
\end{equation}%
Rewriting the left hand side of (\ref{1ad1}),%
\begin{equation*}
\left[ 
\begin{array}{cccccc}
k_{r_{1}}^{2} & 0 & 0 & 0 & 0 & 0 \\ 
0 & k_{\psi _{1}}^{2} & 0 & 0 & 0 & 0 \\ 
0 & 0 & k_{\theta _{1}}^{2} & 0 & 0 & 0 \\ 
0 & 0 & 0 & k_{r_{2}}^{2} & 0 & 0 \\ 
0 & 0 & 0 & 0 & k_{\psi _{2}}^{2} & 0 \\ 
0 & 0 & 0 & 0 & 0 & k_{\theta _{2}}^{2}%
\end{array}%
\right] \left[ 
\begin{array}{c}
\Delta r_{1} \\ 
\Delta \psi _{1} \\ 
\Delta \theta _{1} \\ 
\Delta r_{2} \\ 
\Delta \psi _{2} \\ 
\Delta \theta _{2}%
\end{array}%
\right] =\left[ 
\begin{array}{cc}
D_{1} & 0_{3,3} \\ 
0_{3,3} & D_{2}%
\end{array}%
\right] \left[ 
\begin{array}{c}
\Delta r_{1} \\ 
\Delta \psi _{1} \\ 
\Delta \theta _{1} \\ 
\Delta r_{2} \\ 
\Delta \psi _{2} \\ 
\Delta \theta _{2}%
\end{array}%
\right] 
\end{equation*}%
\begin{equation}
=\left[ \QDATOP{M_{1}}{M_{2}}\right] \left[ 
\begin{array}{c}
a_{1} \\ 
a_{2} \\ 
a_{3}%
\end{array}%
\right] \text{ .}  \label{1aj}
\end{equation}%
Hence, we have%
\begin{equation}
D_{1}\left[ 
\begin{array}{c}
\Delta r_{1} \\ 
\Delta \psi _{1} \\ 
\Delta \theta _{1}%
\end{array}%
\right] =M_{1}\left[ 
\begin{array}{c}
a_{1} \\ 
a_{2} \\ 
a_{3}%
\end{array}%
\right] \text{ ,}  \label{1ak}
\end{equation}%
\begin{equation}
D_{2}\left[ 
\begin{array}{c}
\Delta r_{2} \\ 
\Delta \psi _{2} \\ 
\Delta \theta _{2}%
\end{array}%
\right] =M_{2}\left[ 
\begin{array}{c}
a_{1} \\ 
a_{2} \\ 
a_{3}%
\end{array}%
\right] \text{ ,}  \label{1al}
\end{equation}%
or,%
\begin{equation}
\left[ 
\begin{array}{c}
\Delta r_{2} \\ 
\Delta \psi _{2} \\ 
\Delta \theta _{2}%
\end{array}%
\right] =D_{2}^{-1}M_{2}M_{1}^{-1}D_{1}\left[ 
\begin{array}{c}
\Delta r_{1} \\ 
\Delta \psi _{1} \\ 
\Delta \theta _{1}%
\end{array}%
\right] \text{ .}  \label{1am}
\end{equation}%
Substituting (\ref{1am}) into (\ref{1v}) yields%
\begin{equation}
D_{2}^{-1}M_{2}M_{1}^{-1}D_{1}\left[ 
\begin{array}{c}
\Delta r_{1} \\ 
\Delta \psi _{1} \\ 
\Delta \theta _{1}%
\end{array}%
\right] =A^{-1}\left( p_{T2},\psi _{2},\theta _{2}\right) \left( A\left(
p_{T1},\psi _{1},\theta _{1}\right) \left[ 
\begin{array}{c}
\Delta r_{1} \\ 
\Delta \psi _{1} \\ 
\Delta \theta _{1}%
\end{array}%
\right] +B_{R}\right)   \label{1an}
\end{equation}%
so we get%
\begin{equation}
\left( D_{2}^{-1}M_{2}M_{1}^{-1}D_{1}-A^{-1}\left( p_{T2},\psi _{2},\theta
_{2}\right) A\left( p_{T1},\psi _{1},\theta _{1}\right) \right) \left[ 
\begin{array}{c}
\Delta r_{1} \\ 
\Delta \psi _{1} \\ 
\Delta \theta _{1}%
\end{array}%
\right] =A^{-1}\left( p_{T2},\psi _{2},\theta _{2}\right) B_{R}  \label{1ao}
\end{equation}%
\begin{equation}
\left[ 
\begin{array}{c}
\Delta r_{1} \\ 
\Delta \psi _{1} \\ 
\Delta \theta _{1}%
\end{array}%
\right] =\left( D_{2}^{-1}M_{2}M_{1}^{-1}D_{1}-A^{-1}\left( p_{T2},\psi
_{2},\theta _{2}\right) A\left( p_{T1},\psi _{1},\theta _{1}\right) \right)
^{-1}A^{-1}\left( p_{T2},\psi _{2},\theta _{2}\right) B_{R}\text{ ,}
\label{1ap}
\end{equation}%
which allows us to obtain $\left( \Delta r_{1},\Delta \psi _{1},\Delta
\theta _{1}\right) $. Finally, substituting (\ref{1ap}) into (\ref{1am}) we
get $\left( \Delta r_{2},\Delta \psi _{2},\Delta \theta _{2}\right) $.

\subsection{Numerical Examples}

We illustrate this idea with a few examples.

\subsubsection{a.}

%
%
%
%
%
%
%
%
%
%
%
%
%
%
%
%
%
%
%
%
%
%

\begin{center}
\begin{tabular}{|l|l|}
\hline
INPUT & OUTPUT \\ \hline\hline
$B_{R}$ = [ 200 500 300]' & Cost = 1.6250e+004 \\ 
$p_{T1}$ = 25000 & $\Delta r_{1}$ = -1.7678e+002 \\ 
$\psi _{1}$ = 0 & $\Delta \psi _{1}$ = -1.0000e-002 \\ 
$\theta _{1}$ = 7.8540e-001 & $\Delta \theta _{1}$ = -1.4142e-003 \\ 
$p_{T2}$ = 50000 & $\Delta r_{2}$ = 3.5355e+001 \\ 
$\psi _{2}$ = 0 & $\Delta \psi _{2}$ = 5.0000e-003 \\ 
$\theta _{2}$ = 2.3562e+000 & $\Delta \theta _{2}$ = -3.5355e-003 \\ \hline
\end{tabular}

\begin{tabular}{|l|l|}
\hline
INPUT CONTINUED & OUTPUT\_\_\_\_\_\_\_\_\_\_\_\_\_\_ \\ \hline\hline
\multicolumn{1}{|l|}{$k_{r1}^{2}$ = 2} &  \\ 
\multicolumn{1}{|l|}{$k_{\psi 1}^{2}$ = 1.2500e+009 = $2\ast \mathrm{PT1}%
^{2} $.} &  \\ 
\multicolumn{1}{|l|}{$k_{\theta 1}^{2}$ = 1.2500e+009} &  \\ 
\multicolumn{1}{|l|}{$k_{r2}^{2}$ = 2} &  \\ 
\multicolumn{1}{|l|}{$k_{\psi 2}^{2}$ = 5.0000e+009 = $2\ast \mathrm{PT2}%
^{2} $.} &  \\ 
\multicolumn{1}{|l|}{$k_{\theta 2}^{2}$ = 5.0000e+009} &  \\ \hline
\end{tabular}
\end{center}

\subsubsection{b.}

%
%
%
%
%
%
%
%
%

\begin{center}
\begin{tabular}{|l|l|}
\hline
INPUT Same as with \textbf{a.} but with & OUTPUT \\ \hline\hline
$B_{R}$ = [ 200 0 500]' & Cost = 3.6250e+004 \\ 
& $\Delta r_{1}$ = -2.4749e+002 \\ 
& $\Delta \psi _{1}$ = 0 \\ 
& $\Delta \theta _{1}$ = -4.2426e-003 \\ 
& $\Delta r_{2}$ = 1.0607e+002 \\ 
& $\Delta \psi _{2}$ = 0 \\ 
& $\Delta \theta _{2}$ = -4.9497e-003 \\ \hline
\end{tabular}
\end{center}

\subsubsection{c.}

%
%
%
%
%
%
%
%
%
%

\begin{center}
\begin{tabular}{|l|l|}
\hline
INPUT Same as with \textbf{a.} but with & OUTPUT \\ \hline\hline
\multicolumn{1}{|l|}{$\psi _{2}$ = $\pi $} & Cost = 1.6250e+004 \\ 
\multicolumn{1}{|l|}{$\theta _{2}$ = $\pi /4$} & $\Delta r_{1}$ =
-1.7678e+002 \\ 
\multicolumn{1}{|l|}{} & $\Delta \psi _{1}$ = -1.0000e-002 \\ 
\multicolumn{1}{|l|}{} & $\Delta \theta _{1}$ = -1.4142e-003 \\ 
\multicolumn{1}{|l|}{} & $\Delta r_{2}$ = 3.5355e+001 \\ 
\multicolumn{1}{|l|}{} & $\Delta \psi _{2}$ = -5.0000e-003 \\ 
\multicolumn{1}{|l|}{} & $\Delta \theta _{2}$ = 3.5355e-003 \\ \hline
\end{tabular}
\end{center}

The input for this case is a variation of the input in \textbf{a}. \bigskip
Note that the output is the same as with \textbf{a} except for a sign swap
between $\Delta \psi _{2}$ and $\Delta \theta _{2}$ to account for the
orientation difference of the \textquotedblleft 2\textquotedblright\
coordinates.

\subsubsection{d.}

%
%
%
%
%
%
%
%
%
%

\begin{center}
\begin{tabular}{|l|l|}
\hline
INPUT Same as with \textbf{a.} but with & OUTPUT \\ \hline\hline
$\psi _{2}$ = $\pi /2$ & Cost = 5.5625e+004 \\ 
$\theta _{2}$ = $\pi /4$ & $\Delta r_{1}$ = -1.7678e+002 \\ 
& $\Delta \psi _{1}$ = -1.0000e-002 \\ 
& $\Delta \theta _{1}$ = -1.4142e-003 \\ 
& $\Delta r_{2}$ = 2.8284e+002 \\ 
& $\Delta \psi _{2}$ = -2.0000e-003 \\ 
& $\Delta \theta _{2}$ = -1.4142e-003 \\ \hline
\end{tabular}
\end{center}

\section{An Optimized Reduced-State Filter For Unknown Bias}

A novel technique for calculating a steady-state reduced-order filter to
track a maneuvering target is presented by Mookerjee and Reifler \cite%
{MooRief}. The filter they derive is optimized for performance with a
stochastic acceleration. In this paper, this technique is modified to derive
a steady-state filter that is optimized for performance with a stochastic
measurement bias. Similar to \cite{MooRief}, the filter developed here is a
reduced-state filter. We can see what a reduced-state filter is by
considering \cite{LL/1989} and \cite{LL/1990}. In these reports, we estimate
the position and velocity of an aircraft (a Beechcraft 1900) with DMEs
(distance measuring equipment), an INS (inertial navigation system) and a
barometric altimeter. The filter (in \cite{LL/1989}) and the smoother (in 
\cite{LL/1990}) were designed with a state-to-estimate range bias in each
DME (up to 5 were used), a state-to-estimate INS drift, and a
state-to-estimate bias in the baro. The filter (or smoother) ran with these
additional bias states in tow (i.e., in addition to the position and
velocity states). (The results in \cite{LL/1989} and \cite{LL/1990} achieved
the design goals in position and velocity accuracy.)

It appears likely the design goals of \cite{MooRief} and this paper are
competing design goals. The design methods discussed in \cite{GandS}, which
are based on the Bode gain-phase relationship, possibly could be brought to
bear to quantify a possible trade-off on the design goals of this paper and 
\cite{MooRief}. We don't cover such trade-offs in this paper, but it could
be a problem for future investigations. The classical control concepts of
the sensitivity function and the complimentary sensitivity function come to
mind.

We use discrete time dynamical equations. It is fair to consider our state
and output (dynamical) equations to be the dual (in the control theory
sense) of the state and output equations, Equations \textbf{(8)\footnote{%
In this section, we often refer to equations from \cite{MooRief}. Hence we
adopt the convention that all equation references appearing in bold typeface
are to equations in \cite{MooRief}.}} and \textbf{(5)}. Compared to the
dynamical equations in \cite{MooRief}, we eliminate the unknown acceleration
from the state equation and add an unknown bias in the output (measurement)
equation, the typical dual situation. We have:%

\begin{equation}
x\left( k+1\right) =\Phi x\left( k\right) +m\left( k\right)  \label{2a}
\end{equation}%
\begin{equation}
z\left( k\right) =H\cdot x\left( k\right) +n\left( k\right) +Wu\left(
x\left( k\right) ,\lambda \right) \text{ .}  \label{2b}
\end{equation}%
The state $x(k)$ at time $k$ is of dimension $n$ and the state transition
matrix $\Phi $ is of dimension $n$ by $n$. The output $z(k)$ at time $k$ is
of dimension $q$ and the output matrix $H$ is of dimension $q$ by $n$. The
process noise term $m(k)$ is of dimension $n$ with covariance $Q$. The
measurement noise term $n(k)$ is of dimension $q$ with covariance $N$. The
bias matrix $W$ is $q$ by $m$. The bias function $u$ is $\Re ^{n}\times \Re
^{p}\rightarrow \Re ^{m}$, and we have that the bias $\lambda $ is a $p$%
-dimensional random vector with mean $\overline{\lambda }$ and covariance $%
\Lambda $.

The time update equation, using (\ref{2a}), is simply%
\begin{equation}
\widehat{x}\left( k+1|k\right) =\Phi \widehat{x}\left( k|k\right) \text{ .}
\label{2b_1}
\end{equation}%
The measurement update equation becomes%
\begin{equation}
\widehat{x}\left( k+1|k+1\right) =\widehat{x}\left( k+1|k\right) +K\left(
z\left( k+1\right) -H\widehat{x}\left( k+1|k\right) -Wu\left( \widehat{x}%
\left( k+1|k\right) ,\overline{\lambda }\right) \right) \text{ ,}  \label{2c}
\end{equation}%
where $K$ is the $n$ by $q$ measurement, or Kalman, gain matrix. In the
steady-state case, which is discussed below, the position gain $\alpha $ and
velocity gain $\beta $ substitute for $K$.

\subsection{Filter Development - General Case}

In this subsection, we develop the filter equations for the general case.
The development in this section is (basically) dual (dual in the sense of
control theory) to Section III in \cite{MooRief}. The error is defined as
(we develop the errors analogous to \textbf{(27)} and \textbf{(32)}):%
\begin{equation}
\varepsilon \left( k+1|k+1\right) \equiv x\left( k+1\right) -\widehat{x}%
\left( k+1|k+1\right)  \label{2d}
\end{equation}%
\begin{equation*}
=x\left( k+1\right) -\widehat{x}\left( k+1|k\right) -K\left( z\left(
k+1\right) -H\widehat{x}\left( k+1|k\right) -Wu\left( \widehat{x}\left(
k+1|k\right) ,\overline{\lambda }\right) \right)
\end{equation*}%
\begin{equation*}
=x\left( k+1\right) -\widehat{x}\left( k+1|k\right)
\end{equation*}%
\begin{equation*}
-K\left( Hx\left( k+1\right) +n\left( k+1\right) +Wu\left( x\left(
k+1\right) ,\lambda \right) -H\widehat{x}\left( k+1|k\right) -Wu\left( 
\widehat{x}\left( k+1|k\right) ,\overline{\lambda }\right) \right) \text{ .}
\end{equation*}%
Continuing,%
\begin{equation*}
\varepsilon \left( k+1|k+1\right) =x\left( k+1\right) -KHx\left( k+1\right)
-Kn\left( k+1\right) -KWu\left( x\left( k+1\right) ,\lambda \right)
\end{equation*}%
\begin{equation*}
-\widehat{x}\left( k+1|k\right) +KH\widehat{x}\left( k+1|k\right) +KWu\left( 
\widehat{x}\left( k+1|k\right) ,\overline{\lambda }\right)
\end{equation*}%
\begin{equation*}
=\Phi x\left( k\right) +m\left( k\right) -KH\left( \Phi x\left( k\right)
+m\left( k\right) \right) -Kn\left( k+1\right) -KWu\left( \Phi x\left(
k\right) ,\lambda \right)
\end{equation*}%
\begin{equation*}
-\Phi \widehat{x}\left( k|k\right) +KH\Phi \widehat{x}\left( k|k\right)
+KWu\left( \Phi \widehat{x}\left( k|k\right) ,\overline{\lambda }\right)
\end{equation*}%
\begin{equation*}
=\left( I-KH\right) \Phi \left( x\left( k\right) -\widehat{x}\left(
k|k\right) \right) +\left( I-KH\right) m\left( k\right) -Kn\left( k+1\right)
\end{equation*}%
\begin{equation*}
-KW\left( u\left( \Phi x\left( k\right) ,\lambda \right) -u\left( \Phi 
\widehat{x}\left( k|k\right) ,\overline{\lambda }\right) \right) \text{ .}
\end{equation*}%
So we have%
\begin{equation}
\varepsilon \left( k+1|k+1\right) =L\Phi \varepsilon \left( k|k\right)
+Lm\left( k\right) -K\left( W\Delta u_{k|k}+n\left( k+1\right) \right) \text{
;}  \label{2e}
\end{equation}%
where%
\begin{equation}
L=\left( I-KH\right)  \label{2f}
\end{equation}%
an $n$ by $n$ matrix, and%
\begin{equation}
\Delta u_{k|k}\equiv u\left( \Phi x\left( k\right) ,\lambda \right) -u\left(
\Phi \widehat{x}\left( k|k\right) ,\overline{\lambda }\right) \text{ .}
\label{2g}
\end{equation}%
We can make the linear approximation%
\begin{equation}
\Delta u_{k|k}\approx \left. \frac{\partial u}{\partial x}\right\vert _{x=%
\widehat{x}\left( k|k\right) ,\lambda =\overline{\lambda }}\Phi \Delta
x+\left. \frac{\partial u}{\partial \lambda }\right\vert _{x=\widehat{x}%
\left( k|k\right) ,\lambda =\overline{\lambda }}\Delta \lambda  \label{2h}
\end{equation}%
so%
\begin{equation}
\Delta x=\varepsilon \left( k|k\right) =x\left( k\right) -\widehat{x}\left(
k|k\right) \text{ ,}  \label{2i}
\end{equation}%
and%
\begin{equation}
\Delta \lambda =\lambda -\overline{\lambda }\text{ .}  \label{2j}
\end{equation}%
We obtain the result:%
\begin{equation*}
\varepsilon \left( k+1|k+1\right) =L\Phi \varepsilon \left( k|k\right)
+Lm\left( k\right)
\end{equation*}%
\begin{equation*}
-KW\left( \left. \frac{\partial u}{\partial x}\right\vert _{x=\widehat{x}%
\left( k|k\right) ,\lambda =\overline{\lambda }}\Phi \varepsilon \left(
k|k\right) +\left. \frac{\partial u}{\partial \lambda }\right\vert _{x=%
\widehat{x}\left( k|k\right) ,\lambda =\overline{\lambda }}\Delta \lambda
\right) -Kn\left( k+1\right)
\end{equation*}%
\begin{equation*}
=\left( L-KW\left. \frac{\partial u}{\partial x}\right\vert _{x=\widehat{x}%
\left( k|k\right) ,\lambda =\overline{\lambda }}\right) \Phi \varepsilon
\left( k|k\right) +Lm\left( k\right) -KW\left. \frac{\partial u}{\partial
\lambda }\right\vert _{x=\widehat{x}\left( k|k\right) ,\lambda =\overline{%
\lambda }}\Delta \lambda -Kn\left( k+1\right)
\end{equation*}%
\begin{equation}
=F\varepsilon \left( k|k\right) +Lm\left( k\right) +C\Delta \lambda
-Kn\left( k+1\right)  \label{2k}
\end{equation}%
where%
\begin{equation}
F=\left( L-KW\left. \frac{\partial u}{\partial x}\right\vert _{x=\widehat{x}%
\left( k|k\right) ,\lambda =\overline{\lambda }}\right) \Phi \text{ ,}
\label{2l}
\end{equation}%
an $n$ by $n$ matrix, and%
\begin{equation}
C=-KW\left. \frac{\partial u}{\partial \lambda }\right\vert _{x=\widehat{x}%
\left( k|k\right) ,\lambda =\overline{\lambda }}\text{ ,}  \label{2m}
\end{equation}%
an $n$ by $p$ matrix.

We now implement the observation made in \cite{MooRief} that the error $%
\varepsilon \left( k|k\right) $ may be viewed as consisting of two
components. The first component of error, $\varepsilon ^{\left( 1\right) }$,
is due to the process noise $m$ and the measurement noise $n$. The second
component of error, $\varepsilon ^{\left( 2\right) }$, is due to the
measurement bias. To the extent that the linear approximation is valid, a
linear analysis holds. That is, the two error inputs may be treated in
separate equations by applying the superposition principle of linear
analysis. 
\begin{equation}
\varepsilon ^{\left( 1\right) }\left( k+1|k+1\right) =F\varepsilon ^{\left(
1\right) }\left( k|k\right) +Lm\left( k\right) -Kn\left( k+1\right)
\label{2n}
\end{equation}%
\begin{equation}
\varepsilon ^{\left( 2\right) }\left( k+1|k+1\right) =F\varepsilon ^{\left(
2\right) }\left( k|k\right) +C\cdot \Delta \lambda \text{ .}  \label{2o}
\end{equation}%
These equations are comparable to \textbf{(33)} and \textbf{(34)}.

In addition, we require update equations for the total covariance and the
covariance of $\varepsilon ^{\left( 1\right) }\left( k|k\right) $. Using (%
\ref{2a}), the time update equation for $\varepsilon ^{\left( 1\right)
}\left( k|k\right) $ is%
\begin{equation}
M\left( k+1|k\right) \equiv E\left[ \varepsilon ^{\left( 1\right) }\left(
k+1|k\right) \varepsilon ^{\left( 1\right) }\left( k+1|k\right) ^{\prime }%
\right]  \label{2p}
\end{equation}%
\begin{equation*}
=\Phi M\left( k|k\right) \Phi ^{\prime }+Q\text{ .}
\end{equation*}%
%
%
%
%
%
%
%
%
%
%
%
%
%
%
%
%
%
%
%
%
%
%
%
%
%
%
%
%
%
%
%
%
%
%
%
%
%
%
%
%
%
%
%
%
%
From (\ref{2n}) the combined (measurement and time) update for the
covariance of $\varepsilon ^{\left( 1\right) }\left( k|k\right) $ is 
\begin{equation*}
M\left( k+1|k+1\right) \equiv E\left[ \varepsilon ^{\left( 1\right) }\left(
k+1|k+1\right) \varepsilon ^{\left( 1\right) }\left( k+1|k+1\right) ^{\prime
}\right]
\end{equation*}%
\begin{equation*}
=E\left[ \left( F\varepsilon ^{\left( 1\right) }\left( k|k\right) +Lm\left(
k\right) -Kn\left( k+1\right) \right) \left( F\varepsilon ^{\left( 1\right)
}\left( k|k\right) +Lm\left( k\right) -Kn\left( k+1\right) \right) ^{\prime }%
\right]
\end{equation*}%
and we use $E\left[ n\left( k\right) n\left( l\right) ^{\prime }\right] =0$
for $k\neq l$ giving $E\left[ F\varepsilon ^{\left( 1\right) }\left(
k|k\right) \left( Kn\left( k+1\right) \right) ^{\prime }\right] =0$. Hence,%
\begin{equation*}
M\left( k+1|k+1\right) =E\left[ F\varepsilon ^{\left( 1\right) }\left(
k|k\right) \varepsilon ^{\left( 1\right) }\left( k|k\right) ^{\prime
}F^{\prime }+Lm\left( k\right) m\left( k\right) ^{\prime }L^{\prime
}+Kn\left( k+1\right) n\left( k+1\right) ^{\prime }K^{\prime }\right]
\end{equation*}%
\begin{equation}
=FM\left( k|k\right) F^{\prime }+LQL^{\prime }+KNK^{\prime }\text{ .}
\label{2q}
\end{equation}%
%
%
%
%
%
%
%
%
%
%
%
%
%
%
%
%
%
%
%
%
%
%
%
%
%
%
%
%
%
%
%
%
%
%
%
%
%
%
%
%
%
%
%
%
%
%
Working towards update equations for the total covariance, we firstly define
the $n$ by $p$ matrices $D\left( k|k\right) $ and $D\left( k+1|k\right) $ as%
\begin{equation}
\varepsilon ^{\left( 2\right) }\left( k|k\right) \equiv D\left( k|k\right)
\cdot \Delta \lambda \text{ .}  \label{2r}
\end{equation}%
We can define $D\left( k|k\right) $ in this way since in view of our
linearized analysis, the system output ($\varepsilon ^{\left( 2\right)
}\left( k|k\right) $) is a linear function of the system input ($\Delta
\lambda $). We proceed by defining%
\begin{equation}
D\left( k+1|k\right) \equiv FD\left( k|k\right) \text{ .}  \label{2r_1}
\end{equation}%
In (\ref{2r}), $\varepsilon ^{\left( 2\right) }$ and $\Delta \lambda $ are
known quantities (the equation defines $D\left( k|k\right) $). In (\ref{2r_1}%
), $F$ and $D\left( k|k\right) $ are known quantities. Then, substituting (%
\ref{2r}) into (\ref{2o}), we obtain%
\begin{equation*}
D\left( k+1|k+1\right) \cdot \Delta \lambda =FD\left( k|k\right) \cdot
\Delta \lambda +C\cdot \Delta \lambda
\end{equation*}%
\begin{equation}
=D\left( k+1|k\right) \cdot \Delta \lambda +C\cdot \Delta \lambda \text{ ,}
\label{2s}
\end{equation}%
and subsequently (assuming (\ref{2s}) holds for all $\Delta \lambda $.)%
\begin{equation}
D\left( k+1|k+1\right) =D\left( k+1|k\right) +C\text{ .}  \label{2t}
\end{equation}

Let $S$ be the total error covariance. By superposition, we get the total
error by the addition of the two error terms. We also make the observation
that since the two errors, $\varepsilon ^{\left( 1\right) }$ and $%
\varepsilon ^{\left( 2\right) }$, originate from independent sources they
remain independent for all times $k$. Looking at $S$, 
\begin{equation*}
S\left( k+1|k\right) \equiv E\left[ \varepsilon \left( k+1|k\right)
\varepsilon \left( k+1|k\right) ^{\prime }\right]
\end{equation*}%
\begin{equation*}
=E\left[ \left( \varepsilon ^{\left( 1\right) }\left( k+1|k\right)
+\varepsilon ^{\left( 2\right) }\left( k+1|k\right) \right) \left(
\varepsilon ^{\left( 1\right) }\left( k+1|k\right) +\varepsilon ^{\left(
2\right) }\left( k+1|k\right) \right) ^{\prime }\right]
\end{equation*}%
\begin{equation*}
=E\left[ \varepsilon ^{\left( 1\right) }\left( k+1|k\right) \varepsilon
^{\left( 1\right) }\left( k+1|k\right) ^{\prime }\right] +E\left[
\varepsilon ^{\left( 2\right) }\left( k+1|k\right) \varepsilon ^{\left(
2\right) }\left( k+1|k\right) ^{\prime }\right]
\end{equation*}%
\begin{equation*}
=M\left( k+1|k\right) +E\left[ \varepsilon ^{\left( 2\right) }\left(
k+1|k\right) \varepsilon ^{\left( 2\right) }\left( k+1|k\right) ^{\prime }%
\right]
\end{equation*}%
\begin{equation*}
=M\left( k+1|k\right) +E\left[ \Phi D\left( k|k\right) \Delta \lambda \cdot
\Delta \lambda ^{\prime }D\left( k|k\right) ^{\prime }\Phi ^{\prime }\right] 
\text{ ,}
\end{equation*}%
using (\ref{2a}), (\ref{2b_1}) and (\ref{2r}). Hence,

\begin{equation*}
S\left( k+1|k\right) =M\left( k+1|k\right) +\Phi D\left( k|k\right) E\left[
\Delta \lambda \Delta \lambda ^{\prime }\right] D\left( k|k\right) ^{\prime
}\Phi ^{\prime }\text{ .}
\end{equation*}%
Finally,

\begin{equation}
S\left( k+1|k\right) =M\left( k+1|k\right) +\Phi D\left( k|k\right) \Lambda
D\left( k|k\right) ^{\prime }\Phi ^{\prime }\text{ .}  \label{2u}
\end{equation}%
Basically, this is the same result as \textbf{(19)}.

We next obtain the measurement update for $S$: 
\begin{equation*}
S\left( k+1|k+1\right) \equiv E\left[ \varepsilon \left( k+1|k+1\right)
\varepsilon \left( k+1|k+1\right) ^{\prime }\right]
\end{equation*}%
\begin{equation*}
=E\left[ \left( x\left( k+1\right) -\widehat{x}\left( k+1|k+1\right) \right)
\left( x\left( k+1\right) -\widehat{x}\left( k+1|k+1\right) \right) ^{\prime
}\right]
\end{equation*}%
\begin{equation*}
=E\left[ \left( x\left( k+1\right) -\widehat{x}\left( k+1|k\right) -K\left[
z\left( k+1\right) -H\widehat{x}\left( k+1|k\right) -Wu\left( \widehat{x}%
\left( k+1|k\right) ,\overline{\lambda }\right) \right] \right) \left( 
\mathrm{ditto}\right) ^{\prime }\right]
\end{equation*}%
\begin{equation*}
=E\left[ \left( x\left( k+1\right) -\widehat{x}\left( k+1|k\right) \right.
\right.
\end{equation*}%
\begin{equation*}
\left. \left. -K\left[ Hx\left( k+1\right) +n\left( k+1\right) +Wu\left(
x\left( k+1\right) ,\lambda \right) -H\widehat{x}\left( k+1|k\right)
-Wu\left( \widehat{x}\left( k+1|k\right) ,\overline{\lambda }\right) \right]
\right) \left( \mathrm{ditto}\right) ^{\prime }\right]
\end{equation*}%
\begin{equation*}
=E\left[ \left( x\left( k+1\right) -\widehat{x}\left( k+1|k\right) -KH\left[
x\left( k+1\right) -\widehat{x}\left( k+1|k\right) \right] \right. \right.
\end{equation*}%
\begin{equation*}
\left. \left. -K\left[ n\left( k+1\right) +Wu\left( x\left( k+1\right)
,\lambda \right) -Wu\left( \widehat{x}\left( k+1|k\right) ,\overline{\lambda 
}\right) \right] \right) \left( \mathrm{ditto}\right) ^{\prime }\right]
\end{equation*}%
\begin{equation*}
=E\left[ \left\{ x\left( k+1\right) -\widehat{x}\left( k+1|k\right) -KH\left[
x\left( k+1\right) -\widehat{x}\left( k+1|k\right) \right] \right. \right.
\end{equation*}%
\begin{equation*}
\left. \left. -K\left[ n\left( k+1\right) +W\left( u\left( x\left(
k+1\right) ,\lambda \right) -u\left( \widehat{x}\left( k+1|k\right) ,%
\overline{\lambda }\right) \right) \right] \right\} \left\{ \mathrm{ditto}%
\right\} ^{\prime }\right]
\end{equation*}%
\begin{equation*}
=E\left[ \left\{ x\left( k+1\right) -\widehat{x}\left( k+1|k\right)
-KH\left( x\left( k+1\right) -\widehat{x}\left( k+1|k\right) \right) \right.
\right.
\end{equation*}%
\begin{equation*}
\left. \left. -K\left[ n\left( k+1\right) +W\left( u\left( \Phi x\left(
k\right) ,\lambda \right) -u\left( \Phi \widehat{x}\left( k|k\right) ,%
\overline{\lambda }\right) \right) \right] \right\} \left\{ \mathrm{ditto}%
\right\} ^{\prime }\right]
\end{equation*}%
\begin{equation*}
=E\left[ \left\{ x\left( k+1\right) -\widehat{x}\left( k+1|k\right)
-KH\left( x\left( k+1\right) -\widehat{x}\left( k+1|k\right) \right) -K\left[
n\left( k+1\right) +W\left( \Delta u_{k|k}\right) \right] \right\} \left\{ 
\mathrm{ditto}\right\} ^{\prime }\right] \text{ .}
\end{equation*}%
We take this next step only to the extent of the approximation, 
\begin{equation*}
S\left( k+1|k+1\right) =E\left[ \left\{ \left( I-KH\right) \left( x\left(
k+1\right) -\widehat{x}\left( k+1|k\right) \right) -K\left[ n\left(
k+1\right) +W\left( \frac{\partial u}{\partial x}\Phi \varepsilon \left(
k|k\right) +\frac{\partial u}{\partial \lambda }\Delta \lambda \right) %
\right] \right\} \right.
\end{equation*}%
\begin{equation*}
\times \left. \left\{ \mathrm{ditto}\right\} ^{\prime }\right]
\end{equation*}%
\begin{equation*}
=E\left[ \left\{ \left( I-KH\right) \varepsilon \left( k+1|k\right) -K\left[
n\left( k+1\right) +W\left( \frac{\partial u}{\partial x}\varepsilon \left(
k+1|k\right) +\frac{\partial u}{\partial \lambda }\Delta \lambda \right) %
\right] \right\} \left\{ \mathrm{ditto}\right\} ^{\prime }\right]
\end{equation*}%
\begin{equation*}
=E\left[ \left\{ \left( I-KH-KW\frac{\partial u}{\partial x}\right)
\varepsilon \left( k+1|k\right) -K\left[ n\left( k+1\right) +W\frac{\partial
u}{\partial \lambda }\Delta \lambda \right] \right\} \left\{ \mathrm{ditto}%
\right\} ^{\prime }\right] \text{ .}
\end{equation*}%
Let%
\begin{equation}
\widetilde{H}=H+W\frac{\partial u}{\partial x}  \label{2w}
\end{equation}%
and%
\begin{equation}
\widetilde{N}=N+W\frac{\partial u}{\partial \lambda }\Lambda \frac{\partial u%
}{\partial \lambda }^{\prime }W^{\prime }\text{ .}  \label{2x}
\end{equation}%
Then%
\begin{equation*}
S\left( k+1|k+1\right) =\left( I-K\widetilde{H}\right) S\left( k+1|k\right)
\left( I-K\widetilde{H}\right) ^{\prime }+K\widetilde{N}K^{\prime }
\end{equation*}%
\begin{equation*}
-\left( I-K\widetilde{H}\right) E\left[ \varepsilon \left( k+1|k\right)
\Delta \lambda ^{\prime }\right] \left( W\frac{\partial u}{\partial \lambda }%
\right) ^{\prime }K^{\prime }
\end{equation*}%
\begin{equation*}
-K\left( W\frac{\partial u}{\partial \lambda }\right) E\left[ \Delta \lambda
\varepsilon \left( k+1|k\right) ^{\prime }\right] \left( I-K\widetilde{H}%
\right) ^{\prime }\text{ .}
\end{equation*}%
Now%
\begin{equation*}
E\left[ \varepsilon \left( k+1|k\right) \Delta \lambda ^{\prime }\right] =E%
\left[ \left( \varepsilon ^{(1)}\left( k+1|k\right) +\varepsilon
^{(2)}\left( k+1|k\right) \right) \Delta \lambda ^{\prime }\right]
\end{equation*}%
\begin{equation*}
=E\left[ \varepsilon ^{(2)}\left( k+1|k\right) \Delta \lambda ^{\prime }%
\right] =\Phi E\left[ \varepsilon ^{(2)}\left( k|k\right) \Delta \lambda
^{\prime }\right]
\end{equation*}%
\begin{equation*}
=\Phi E\left[ D\left( k|k\right) \Delta \lambda \Delta \lambda ^{\prime }%
\right] =\Phi D\left( k|k\right) \Lambda \text{ .}
\end{equation*}%
Hence,%
\begin{equation*}
S\left( k+1|k+1\right) =\left( I-K\widetilde{H}\right) S\left( k+1|k\right)
\left( I-K\widetilde{H}\right) ^{\prime }+K\widetilde{N}K^{\prime }
\end{equation*}%
\begin{equation}
-\left( I-K\widetilde{H}\right) \Phi D\left( k|k\right) \Lambda \left( W%
\frac{\partial u}{\partial \lambda }\right) ^{\prime }K^{\prime }-K\left( W%
\frac{\partial u}{\partial \lambda }\right) \Lambda D\left( k|k\right)
^{\prime }\Phi ^{\prime }\left( I-K\widetilde{H}\right) ^{\prime }\text{ .}
\label{2v}
\end{equation}%
Equation (\ref{2v}) is similar in form as \textbf{(37)}. We select $K$ so as
to minimize the trace of $S\left( k+1|k+1\right) $: $tr\left( S\left(
k+1|k+1\right) \right) $. We minimize $tr\left( S\left( k+1|k+1\right)
\right) $ because for positive definite matrices trace going to zero implies
the matrix $L_{2}$ norm goes to zero. Let $P$ be a positive definite $n$ by $%
n$ matrix, we have%
\begin{equation*}
\left\Vert P\right\Vert \leq tr(P)\leq n\,\left\Vert P\right\Vert \text{ .}
\end{equation*}%
So minimizing $tr\left( P\right) $, gives a smaller upper bound for $%
\left\Vert P\right\Vert $. We find the optimal $K$ ($K$ that minimizes $%
tr\left( S\left( k+1|k+1\right) \right) $) by taking derivatives with
respect to $K$, setting the result to zero and solving for $K$. We recall
the following facts: Let $A$ be a matrix independent of $K$,%
\begin{equation*}
\frac{\partial }{\partial K}tr\left( KAK^{\prime }\right) =K\left(
A+A^{\prime }\right) \text{ ;}
\end{equation*}%
\begin{equation*}
\frac{\partial }{\partial K}tr\left( KA\right) =A^{\prime }\text{ ;}
\end{equation*}%
\begin{equation*}
\frac{\partial }{\partial K}tr\left( K^{\prime }A\right) =A\text{ ;}
\end{equation*}%
\begin{equation*}
\frac{\partial }{\partial K}tr\left( AK\right) =A^{\prime }\text{ ;}
\end{equation*}%
\begin{equation*}
\frac{\partial }{\partial K}tr\left( AK^{\prime }\right) =A\text{ .}
\end{equation*}%
Using these facts on the terms of (\ref{2v}), 
\begin{equation*}
\frac{\partial }{\partial K}tr\left( K\widetilde{H}S\left( k+1|k\right) 
\widetilde{H}^{\prime }K^{\prime }\right) =2K\widetilde{H}S\left(
k+1|k\right) \widetilde{H}^{\prime }\text{ ;}
\end{equation*}%
\begin{equation*}
\frac{\partial }{\partial K}tr\left( K\widetilde{N}K^{\prime }\right) =2K%
\widetilde{N}\text{ ;}
\end{equation*}%
\begin{equation*}
\frac{\partial }{\partial K}tr\left( S\left( k+1|k\right) \right) =0\text{ ;}
\end{equation*}%
\begin{equation*}
\frac{\partial }{\partial K}tr\left( S\left( k+1|k\right) \widetilde{H}%
^{\prime }K^{\prime }\right) =S\left( k+1|k\right) \widetilde{H}^{\prime }%
\text{ ;}
\end{equation*}%
\begin{equation*}
\frac{\partial }{\partial K}tr\left( K\widetilde{H}S\left( k+1|k\right)
\right) =S\left( k+1|k\right) \widetilde{H}^{\prime }\text{ ;}
\end{equation*}%
\begin{equation*}
\frac{\partial }{\partial K}tr\left( \Phi D\left( k|k\right) \Lambda \left( W%
\frac{\partial u}{\partial \lambda }\right) ^{\prime }K^{\prime }\right)
=\Phi D\left( k|k\right) \Lambda \left( W\frac{\partial u}{\partial \lambda }%
\right) ^{\prime }\text{ ;}
\end{equation*}%
\begin{equation*}
\frac{\partial }{\partial K}tr\left( K\widetilde{H}\Phi D\left( k|k\right)
\Lambda \left( W\frac{\partial u}{\partial \lambda }\right) ^{\prime
}K^{\prime }\right) =K\left( \widetilde{H}\Phi D\left( k|k\right) \Lambda
\left( W\frac{\partial u}{\partial \lambda }\right) ^{\prime }+\left( W\frac{%
\partial u}{\partial \lambda }\right) \Lambda D\left( k|k\right) ^{\prime
}\Phi ^{\prime }\widetilde{H}^{\prime }\right) \text{ ;}
\end{equation*}%
\begin{equation*}
\frac{\partial }{\partial K}tr\left( K\left( W\frac{\partial u}{\partial
\lambda }\right) \Lambda D\left( k|k\right) ^{\prime }\Phi ^{\prime }\right)
=\Phi D\left( k|k\right) \Lambda \left( W\frac{\partial u}{\partial \lambda }%
\right) ^{\prime }\text{ ;}
\end{equation*}%
\begin{equation*}
\frac{\partial }{\partial K}tr\left( K\left( W\frac{\partial u}{\partial
\lambda }\right) \Lambda D\left( k|k\right) ^{\prime }\Phi ^{\prime }%
\widetilde{H}^{\prime }K^{\prime }\right) =K\left( \left( W\frac{\partial u}{%
\partial \lambda }\right) \Lambda D\left( k|k\right) ^{\prime }\Phi ^{\prime
}\widetilde{H}^{\prime }+\widetilde{H}\Phi D\left( k|k\right) \Lambda \left(
W\frac{\partial u}{\partial \lambda }\right) ^{\prime }\right) \text{ .}
\end{equation*}%
We differentiate the trace of $S\left( k+1|k+1\right) $ as it is represented
in (\ref{2v}) by $K$ and set the result equal to zero,%
\begin{equation*}
\frac{\partial }{\partial K}tr\left( S\left( k+1|k+1\right) \right) =2K%
\widetilde{H}S\left( k+1|k\right) \widetilde{H}^{\prime }+2K\widetilde{N}%
+0-S\left( k+1|k\right) \widetilde{H}^{\prime }-S\left( k+1|k\right) 
\widetilde{H}^{\prime }
\end{equation*}%
\begin{equation*}
-\Phi D\left( k|k\right) \Lambda \left( W\frac{\partial u}{\partial \lambda }%
\right) ^{\prime }+K\left( \widetilde{H}\Phi D\left( k|k\right) \Lambda
\left( W\frac{\partial u}{\partial \lambda }\right) ^{\prime }+\left( W\frac{%
\partial u}{\partial \lambda }\right) \Lambda D\left( k|k\right) ^{\prime
}\Phi ^{\prime }\widetilde{H}^{\prime }\right)
\end{equation*}%
\begin{equation*}
-\Phi D\left( k|k\right) \Lambda \left( W\frac{\partial u}{\partial \lambda }%
\right) ^{\prime }+K\left( \left( W\frac{\partial u}{\partial \lambda }%
\right) \Lambda D\left( k|k\right) ^{\prime }\Phi ^{\prime }\widetilde{H}%
^{\prime }+\widetilde{H}\Phi D\left( k|k\right) \Lambda \left( W\frac{%
\partial u}{\partial \lambda }\right) ^{\prime }\right) =0\text{ .}
\end{equation*}%
After some algebra, 
\begin{equation*}
K\left( \widetilde{H}S\left( k+1|k\right) \widetilde{H}^{\prime }+\widetilde{%
N}+\widetilde{H}\Phi D\left( k|k\right) \Lambda \left( W\frac{\partial u}{%
\partial \lambda }\right) ^{\prime }+\left( W\frac{\partial u}{\partial
\lambda }\right) \Lambda D\left( k|k\right) ^{\prime }\Phi ^{\prime }%
\widetilde{H}^{\prime }\right)
\end{equation*}%
\begin{equation}
=S\left( k+1|k\right) \widetilde{H}^{\prime }+\Phi D\left( k|k\right)
\Lambda \left( W\frac{\partial u}{\partial \lambda }\right) ^{\prime }\text{
.}  \label{2x_1}
\end{equation}%
The optimal filter gain is%
\begin{equation*}
K=\left( \widetilde{H}S\left( k+1|k\right) \widetilde{H}^{\prime }+%
\widetilde{N}+\widetilde{H}\Phi D\left( k|k\right) \Lambda \left( W\frac{%
\partial u}{\partial \lambda }\right) ^{\prime }+\left( W\frac{\partial u}{%
\partial \lambda }\right) \Lambda D\left( k|k\right) ^{\prime }\Phi ^{\prime
}\widetilde{H}^{\prime }\right) ^{-1}
\end{equation*}%
\begin{equation}
\times \left( S\left( k+1|k\right) \widetilde{H}^{\prime }+\Phi D\left(
k|k\right) \Lambda \left( W\frac{\partial u}{\partial \lambda }\right)
^{\prime }\right) \text{ ,}  \label{2y}
\end{equation}%
which is an $n$ by $q$ matrix.

\subsection{Filter Development - Steady-State Case}

We next examine the steady-state case of our problem. Referencing equations 
\textbf{(41-43)} we have%
\begin{equation}
\Phi =\left[ 
\begin{array}{cc}
1 & T \\ 
0 & 1%
\end{array}%
\right]  \label{2z}
\end{equation}%
%
%
%
%
%
%
%
%
%
%
%
%
%
%
%
%
%
%
%
%
%
%
%
%
%
%
%
%
%
%
%
%
%
%
%
%
%
\begin{equation}
H=\left[ 
\begin{array}{cc}
1 & 0%
\end{array}%
\right]  \label{2ab}
\end{equation}%
\begin{equation}
W=1\text{ .}  \label{2ac}
\end{equation}%
Hence, with $p$ as position, $v$ as velocity, and $z$ as the measurement,
the state transition and output equations for the steady-state case are 
\begin{equation}
\left[ 
\begin{array}{c}
p\left( k+1\right) \\ 
v\left( k+1\right)%
\end{array}%
\right] =\left[ 
\begin{array}{cc}
1 & T \\ 
0 & 1%
\end{array}%
\right] \left[ 
\begin{array}{c}
p\left( k\right) \\ 
v\left( k\right)%
\end{array}%
\right] +\left[ 
\begin{array}{c}
0 \\ 
1%
\end{array}%
\right] m\left( k\right)  \label{2ad}
\end{equation}%
\begin{equation}
z\left( k\right) =\left[ 
\begin{array}{cc}
1 & 0%
\end{array}%
\right] \left[ 
\begin{array}{c}
p\left( k\right) \\ 
v\left( k\right)%
\end{array}%
\right] +n\left( k\right) +u\left( x\left( k\right) ,\lambda \left( k\right)
\right) \text{ .}  \label{2ae}
\end{equation}%
Setting $u\left( x,\lambda \right) =\lambda $, our linear approximations
from (\ref{2h}) become%
\begin{equation}
\left. \frac{\partial u}{\partial x}\right\vert _{x=\widehat{x}\left(
k|k\right) ,\lambda =\overline{\lambda }}=\left[ 
\begin{array}{cc}
0 & 0%
\end{array}%
\right]  \label{2af}
\end{equation}%
and
\begin{equation}
\left. \frac{\partial u}{\partial \lambda }\right\vert _{x=\widehat{x}\left(
k|k\right) ,\lambda =\overline{\lambda }}=1\text{ .}  \label{2ag}
\end{equation}%
Then, substituting (\ref{2af}) into (\ref{2w})%
\begin{equation}
\widetilde{H}=H+W\frac{\partial u}{\partial x}=\left[ 
\begin{array}{cc}
1 & 0%
\end{array}%
\right] +1\cdot \left[ 
\begin{array}{cc}
0 & 0%
\end{array}%
\right] =\left[ 
\begin{array}{cc}
1 & 0%
\end{array}%
\right] =H  \label{2ah}
\end{equation}%
and substituting (\ref{2ag}) into (\ref{2x})%
\begin{equation}
\widetilde{N}=N+W\frac{\partial u}{\partial \lambda }\Lambda \frac{\partial u%
}{\partial \lambda }^{\prime }W^{\prime }=N+1\cdot 1\cdot \Lambda \cdot
1\cdot 1=N+\Lambda \text{ .}  \label{2ai}
\end{equation}%
We have that the steady-state filter gain is
\begin{equation}
\overline{K}\equiv \left[ 
\begin{array}{c}
\alpha \\ 
\beta /T%
\end{array}%
\right] \text{ .}  \label{2aj}
\end{equation}%
As mentioned previously, (\ref{2aj}) is where $\alpha $ and $\beta $ fit in
for the Kalman gain matrix $K$ as given by (\ref{2y}). These gains are
obtained by computing the steady-state values for all the variables in (\ref%
{2y}). The major objective of this section is to find a relationship between 
$\alpha $ and $\beta $.

The steady-state version of $L$ from (\ref{2f}), $\overline{L}$, is%
\begin{equation}
\overline{L}=\left( I-\overline{K}H\right) =\left[ 
\begin{array}{cc}
1 & 0 \\ 
0 & 1%
\end{array}%
\right] -\left[ 
\begin{array}{c}
\alpha \\ 
\beta /T%
\end{array}%
\right] \left[ 
\begin{array}{cc}
1 & 0%
\end{array}%
\right] =\left[ 
\begin{array}{cc}
1-\alpha & 0 \\ 
-\beta /T & 1%
\end{array}%
\right] \text{ .}  \label{2ak}
\end{equation}%
The steady-state version of $F$ from (\ref{2l}), $\overline{F}$, is%
\begin{equation}
\overline{F}=\left( \left[ 
\begin{array}{cc}
1-\alpha & 0 \\ 
-\beta /T & 1%
\end{array}%
\right] -\left[ 
\begin{array}{c}
\alpha \\ 
\beta /T%
\end{array}%
\right] \cdot 1\cdot \left[ 
\begin{array}{cc}
0 & 0%
\end{array}%
\right] \right) \left[ 
\begin{array}{cc}
1 & T \\ 
0 & 1%
\end{array}%
\right] =\left[ 
\begin{array}{cc}
1-\alpha & \left( 1-\alpha \right) T \\ 
-\beta /T & 1-\beta%
\end{array}%
\right] \text{ .}  \label{2al}
\end{equation}%
The eigenvalues of $\overline{F}$ are 
\begin{equation}
\lambda _{1,2}=1-\frac{\left( \alpha +\beta \right) }{2}\pm \frac{1}{2}\sqrt{%
2\alpha \beta -4\beta +\alpha ^{2}+\beta ^{2}}\text{ .}  \label{2al1}
\end{equation}%
Then, referring to (\ref{2m}), the steady-state version of $C$, $\overline{C}
$, is%
\begin{equation}
\overline{C}=-\left[ 
\begin{array}{c}
\alpha \\ 
\beta /T%
\end{array}%
\right] \cdot 1\cdot 1=-\left[ 
\begin{array}{c}
\alpha \\ 
\beta /T%
\end{array}%
\right] \text{ .}  \label{2am}
\end{equation}%
%
%
%
%
%
%
%
%
%
%
%
%
%
%
%
%
%
%
%
%
%
%
%
%
%
%
%
%
%
%
%
%
%
%
%
%
%
%
%
%
%
%
%
%
%
%
%
The measurement updated steady-state covariance $M$, referring to (\ref{2q}%
), is 
\begin{equation}
\overline{M}\equiv \lim_{k\rightarrow \infty }M\left( k|k\right) =\overline{F%
}\,\overline{M}\,\overline{F}^{\prime }+\overline{L}Q\overline{L}^{\prime }+%
\overline{K}\,N\,\overline{K}^{\prime }  \label{2am1}
\end{equation}%
Superimpose the two noise terms by letting 
\begin{equation}
\overline{M}=\overline{M}_{Q}+\overline{M}_{N}  \label{2am1o}
\end{equation}%
and solve:%
\begin{equation}
\overline{M}_{N}=\overline{F}\,\overline{M}_{N}\,\overline{F}^{\prime }+%
\overline{K}\,N\,\overline{K}^{\prime }  \label{2am1i}
\end{equation}%
\begin{equation}
\overline{M}_{Q}=\overline{F}\,\overline{M}_{Q}\,\overline{F}^{\prime }+%
\overline{L}Q\overline{L}^{\prime }\text{ .}  \label{2am1ii}
\end{equation}%
Comparing (\ref{2aj}), (\ref{2al}) and (\ref{2am1i}) to \textbf{(46)}, 
\textbf{(47)} and \textbf{(48)}, we see that our solution for $\overline{M}%
_{N}$ is of the same form as \textbf{(49).} Consequently, 
\begin{equation}
\overline{M}_{N}=\frac{N}{\alpha \left( 4-2\alpha -\beta \right) }\left[ 
\begin{array}{cc}
2\alpha ^{2}+2\beta -3\alpha \beta & \beta \left( 2\alpha -\beta \right) /T
\\ 
\beta \left( 2\alpha -\beta \right) /T & 2\beta ^{2}/T^{2}%
\end{array}%
\right] \text{ .}  \label{2am1iii}
\end{equation}%
The solution of $\overline{M}_{Q}$ remains to be determined. Substituting (%
\ref{2ak}), (\ref{2al}) and (\ref{2ad}) into (\ref{2am1ii}) gives%
\begin{equation*}
\overline{M}_{Q}=\left[ 
\begin{array}{cc}
1-\alpha & \left( 1-\alpha \right) T \\ 
-\beta /T & 1-\beta%
\end{array}%
\right] \overline{M}_{Q}\left[ 
\begin{array}{cc}
1-\alpha & -\beta /T \\ 
\left( 1-\alpha \right) T & 1-\beta%
\end{array}%
\right] +\left[ 
\begin{array}{cc}
1-\alpha & 0 \\ 
-\beta /T & 1%
\end{array}%
\right] \left[ 
\begin{array}{cc}
0 & 0 \\ 
0 & q_{22}%
\end{array}%
\right] \left[ 
\begin{array}{cc}
1-\alpha & -\beta /T \\ 
0 & 1%
\end{array}%
\right] \text{ .}
\end{equation*}%
In longhand,%
\begin{equation*}
\overline{M}_{Q}=\left[ 
\begin{array}{cc}
\overline{m}_{11Q} & \overline{m}_{12Q} \\ 
\overline{m}_{12Q} & \overline{m}_{22Q}%
\end{array}%
\right]
\end{equation*}%
\begin{equation*}
=\left[ 
\begin{array}{c}
\left( 1-\alpha \right) ^{2}\left( \overline{m}_{11Q}+2T\overline{m}%
_{12Q}+T^{2}\overline{m}_{22Q}\right) \\ 
\left( 1-\alpha \right) \left( -\beta \overline{m}_{11Q}/T+\left( 1-2\beta
\right) \overline{m}_{12Q}+T\left( 1-\beta \right) \overline{m}_{22Q}\right)%
\end{array}%
\right.
\end{equation*}%
\begin{equation}
\left. 
\begin{array}{c}
\left( 1-\alpha \right) \left( -\beta \overline{m}_{11Q}/T+\left( 1-2\beta
\right) \overline{m}_{12Q}+T\left( 1-\beta \right) \overline{m}_{22Q}\right)
\\ 
\left( \beta ^{2}/T^{2}\right) \overline{m}_{11Q}+\left( 2\beta \left( \beta
-1\right) /T\right) \overline{m}_{12Q}+\left( 1-2\beta +\beta ^{2}\right) 
\overline{m}_{22Q}%
\end{array}%
\right] +\allowbreak \left[ 
\begin{array}{cc}
0 & 0 \\ 
0 & q_{22}%
\end{array}%
\right] \text{ .}  \label{2am2}
\end{equation}%
Hence,%
\begin{equation*}
\left[ 
\begin{array}{c}
\overline{m}_{11Q}-\left( 1-\alpha \right) ^{2}\left( \overline{m}_{11Q}+2T%
\overline{m}_{12Q}+T^{2}\overline{m}_{22Q}\right) \\ 
\overline{m}_{12Q}-\left( 1-\alpha \right) \left( -\beta \overline{m}%
_{11Q}/T+\left( 1-2\beta \right) \overline{m}_{12Q}+T\left( 1-\beta \right) 
\overline{m}_{22Q}\right)%
\end{array}%
\right.
\end{equation*}%
\begin{equation*}
\left. 
\begin{array}{c}
\overline{m}_{12Q}-\left( 1-\alpha \right) \left( -\beta \overline{m}%
_{11Q}/T+\left( 1-2\beta \right) \overline{m}_{12Q}+T\left( 1-\beta \right) 
\overline{m}_{22Q}\right) \\ 
\overline{m}_{22Q}-\left( \beta ^{2}/T^{2}\right) \overline{m}_{11Q}-\left(
2\beta \left( \beta -1\right) /T\right) \overline{m}_{12Q}-\left( 1-2\beta
+\beta ^{2}\right) \overline{m}_{22Q}%
\end{array}%
\right]
\end{equation*}%
\begin{equation}
=\left[ 
\begin{array}{cc}
0 & 0 \\ 
0 & q_{22}%
\end{array}%
\right] \text{ .}  \label{2am3}
\end{equation}%
We get three equations in the three unknowns $m_{11Q}$, $m_{12Q}$ and $%
m_{22Q}$,%
\begin{equation*}
\left[ 
\begin{array}{ccc}
1-\left( 1-\alpha \right) ^{2} & -2\left( 1-\alpha \right) ^{2}T & -\left(
1-\alpha \right) ^{2}T^{2} \\ 
\left( 1-\alpha \right) \beta /T & 1-\left( 1-\alpha \right) \left( 1-2\beta
\right) & -\left( 1-\alpha \right) \left( 1-\beta \right) T \\ 
-\left( \beta ^{2}/T^{2}\right) & -2\beta \left( \beta -1\right) /T & 
1-\left( 1-2\beta +\beta ^{2}\right)%
\end{array}%
\right] \left[ 
\begin{array}{c}
\overline{m}_{11Q} \\ 
\overline{m}_{12Q} \\ 
\overline{m}_{22Q}%
\end{array}%
\right] \text{ }
\end{equation*}%
\begin{equation*}
=\left[ 
\begin{array}{c}
0 \\ 
0 \\ 
q_{22}%
\end{array}%
\right] \text{ .}
\end{equation*}%
Taking the matrix inverse to solve for $\overline{M}_{Q}$,%
\begin{equation*}
\left[ 
\begin{array}{c}
\overline{m}_{11Q} \\ 
\overline{m}_{12Q} \\ 
\overline{m}_{22Q}%
\end{array}%
\right] =\left[ 
\begin{array}{ccc}
1-\left( 1-\alpha \right) ^{2} & -2\left( 1-\alpha \right) ^{2}T & -\left(
1-\alpha \right) ^{2}T^{2} \\ 
\left( 1-\alpha \right) \beta /T & 1-\left( 1-\alpha \right) \left( 1-2\beta
\right) & -\left( 1-\alpha \right) \left( 1-\beta \right) T \\ 
-\left( \beta ^{2}/T^{2}\right) & -2\beta \left( \beta -1\right) /T & 
1-\left( 1-2\beta +\beta ^{2}\right)%
\end{array}%
\right] ^{-1}\left[ 
\begin{array}{c}
0 \\ 
0 \\ 
q_{22}%
\end{array}%
\right] \text{ .}
\end{equation*}%
The determinant of this matrix, $4\alpha \beta -\alpha \beta ^{2}-2\alpha
^{2}\beta =\alpha \beta \left( 4-\beta -2\alpha \right) $, should not be
zero for the inverse to exist. This is satisfied by these conditions:%
\begin{equation}
\begin{array}{ll}
\mathrm{1.} & \alpha \neq 0 \\ 
\mathrm{2.} & \beta \neq 0 \\ 
\mathrm{3.} & \beta \neq 4-2\alpha%
\end{array}%
\text{ .}  \label{2am4}
\end{equation}%
If the determinant is not zero, we can obtain the solution:%
\begin{equation}
\left[ 
\begin{array}{c}
\overline{m}_{11Q} \\ 
\overline{m}_{12Q} \\ 
\overline{m}_{22Q}%
\end{array}%
\right] =q_{22}\cdot \left[ 
\begin{array}{c}
T^{2}\left( -2+5\alpha -4\alpha ^{2}+\alpha ^{3}\right) \\ 
T\left( -2\alpha +\beta -\alpha \beta +3\alpha ^{2}-\alpha ^{3}\right) \\ 
\left( -2\beta +2\alpha \beta -2\alpha ^{2}+\alpha ^{3}\right)%
\end{array}%
\right] /\left( -4\alpha \beta +\alpha \beta ^{2}+2\alpha ^{2}\beta \right) 
\text{ .}  \label{2an}
\end{equation}%
In matrix form,%
\begin{equation}
\overline{M}_{Q}=\frac{q_{22}}{\left( -4\alpha \beta +\alpha \beta
^{2}+2\alpha ^{2}\beta \right) }\left[ 
\begin{array}{cc}
T^{2}\left( -2+5\alpha -4\alpha ^{2}+\alpha ^{3}\right) & T\left( -2\alpha
+\beta -\alpha \beta +3\alpha ^{2}-\alpha ^{3}\right) \\ 
T\left( -2\alpha +\beta -\alpha \beta +3\alpha ^{2}-\alpha ^{3}\right) & 
\left( -2\beta +2\alpha \beta -2\alpha ^{2}+\alpha ^{3}\right)%
\end{array}%
\right] \text{ .}  \label{2an1}
\end{equation}%
And finally $\overline{M}$ is obtained from (\ref{2am1o}), (\ref{2am1iii})
and (\ref{2an1}):%
\begin{equation*}
\overline{M}=\frac{N}{\alpha \left( 4-2\alpha -\beta \right) }\left[ 
\begin{array}{cc}
2\alpha ^{2}+2\beta -3\alpha \beta & \beta \left( 2\alpha -\beta \right) /T
\\ 
\beta \left( 2\alpha -\beta \right) /T & 2\beta ^{2}/T^{2}%
\end{array}%
\right]
\end{equation*}%
\begin{equation}
+\frac{q_{22}}{\left( -4\alpha \beta +\alpha \beta ^{2}+2\alpha ^{2}\beta
\right) }\left[ 
\begin{array}{cc}
T^{2}\left( -2+5\alpha -4\alpha ^{2}+\alpha ^{3}\right) & T\left( -2\alpha
+\beta -\alpha \beta +3\alpha ^{2}-\alpha ^{3}\right) \\ 
T\left( -2\alpha +\beta -\alpha \beta +3\alpha ^{2}-\alpha ^{3}\right) & 
\left( -2\beta +2\alpha \beta -2\alpha ^{2}+\alpha ^{3}\right)%
\end{array}%
\right] \text{ .}  \label{2an2}
\end{equation}

%
%
%
%
%
%
%
%
%
%
%
%
%
%
%
%
%
%
%
%
%
%
%
We see that $\overline{m}_{11}=\overline{m}_{11}\left( \alpha ,\beta
,T,N,q_{22}\right) $, $\overline{m}_{12}=\overline{m}_{12}\left( \alpha
,\beta ,T,N,q_{22}\right) $ and $\overline{m}_{22}=\overline{m}_{22}\left(
\alpha ,\beta ,T,N,q_{22}\right) $. %
%
%
%
The usual technique for solving the Liapunov equation (\ref{2am1}) is by
algebraic manipulation and using the symmetry of the matrix, as demonstrated
with the solution (\ref{2an}). Numerical solutions may be obtained by
repeated propagation until steady-state is arrived at.

The time updated steady-state covariance $M$, referring to (\ref{2p}), is 
\begin{equation}
\overset{\cdot }{M}\equiv \lim_{k\rightarrow \infty }M\left( k+1|k\right)
=\lim_{k\rightarrow \infty }\Phi M\left( k|k\right) \Phi ^{\prime }+Q=\Phi
\left( \overline{M}_{N}+\overline{M}_{Q}\right) \Phi ^{\prime }+\allowbreak Q%
\text{ .}  \label{2ao}
\end{equation}%
We get%
\begin{equation*}
\Phi \overline{M}_{N}\Phi ^{\prime }=\frac{N}{\alpha \left( 4-2\alpha -\beta
\right) }\left[ 
\begin{array}{cc}
1 & T \\ 
0 & 1%
\end{array}%
\right] \left[ 
\begin{array}{cc}
2\alpha ^{2}+2\beta -3\alpha \beta & \beta \left( 2\alpha -\beta \right) /T
\\ 
\beta \left( 2\alpha -\beta \right) /T & 2\beta ^{2}/T^{2}%
\end{array}%
\right] \left[ 
\begin{array}{cc}
1 & 0 \\ 
T & 1%
\end{array}%
\right]
\end{equation*}%
\begin{equation*}
=\frac{N}{\alpha \left( 4-2\alpha -\beta \right) }\left[ 
\begin{array}{cc}
2\alpha ^{2}+2\beta +\alpha \beta & \beta \left( 2\alpha +\beta \right) /T
\\ 
\beta \left( 2\alpha +\beta \right) /T & 2\beta ^{2}/T^{2}%
\end{array}%
\right]
\end{equation*}%
and%
\begin{equation*}
\Phi \overline{M}_{Q}\Phi ^{\prime }=
\end{equation*}%
\begin{equation*}
=\frac{q_{22}}{\left( -4\alpha \beta +\alpha \beta ^{2}+2\alpha ^{2}\beta
\right) }\left[ 
\begin{array}{cc}
1 & T \\ 
0 & 1%
\end{array}%
\right] \left[ 
\begin{array}{cc}
T^{2}\left( -2+5\alpha -4\alpha ^{2}+\alpha ^{3}\right) & T\left( -2\alpha
+\beta -\alpha \beta +3\alpha ^{2}-\alpha ^{3}\right) \\ 
T\left( -2\alpha +\beta -\alpha \beta +3\alpha ^{2}-\alpha ^{3}\right) & 
\left( -2\beta +2\alpha \beta -2\alpha ^{2}+\alpha ^{3}\right)%
\end{array}%
\right] \left[ 
\begin{array}{cc}
1 & 0 \\ 
T & 1%
\end{array}%
\right]
\end{equation*}%
\begin{equation*}
=\frac{q_{22}}{\left( -4\alpha \beta +\alpha \beta ^{2}+2\alpha ^{2}\beta
\right) }\left[ 
\begin{array}{cc}
T^{2}\left( -2+\alpha \right) & T\left( -2\alpha -\beta +\alpha \beta
+\alpha ^{2}\right) \\ 
T\left( -2\alpha -\beta +\alpha \beta +\alpha ^{2}\right) & \left( -2\beta
+2\alpha \beta -2\alpha ^{2}+\alpha ^{3}\right)%
\end{array}%
\right] \allowbreak \text{ .}
\end{equation*}%
Hence,%
\begin{equation*}
\overset{\cdot }{M}=\frac{N}{\alpha \left( 4-2\alpha -\beta \right) }\left[ 
\begin{array}{cc}
2\alpha ^{2}+2\beta +\alpha \beta & \beta \left( 2\alpha +\beta \right) /T
\\ 
\beta \left( 2\alpha +\beta \right) /T & 2\beta ^{2}/T^{2}%
\end{array}%
\right]
\end{equation*}%
\begin{equation}
+\frac{q_{22}}{\left( -4\alpha \beta +\alpha \beta ^{2}+2\alpha ^{2}\beta
\right) }\left[ 
\begin{array}{cc}
T^{2}\left( -2+\alpha \right) & T\left( -2\alpha -\beta +\alpha \beta
+\alpha ^{2}\right) \\ 
T\left( -2\alpha -\beta +\alpha \beta +\alpha ^{2}\right) & \left( -2\beta
+2\alpha \beta -2\alpha ^{2}+\alpha ^{3}\right)%
\end{array}%
\right] +\left[ 
\begin{array}{cc}
0 & 0 \\ 
0 & q_{22}%
\end{array}%
\right] \text{ .}  \label{2ao1}
\end{equation}

We need steady-state versions of these: designating%
\begin{equation*}
\overline{D}\equiv \lim_{k\rightarrow \infty }D\left( k|k\right)
\end{equation*}%
\begin{equation*}
\overset{\cdot }{D}\equiv \lim_{k\rightarrow \infty }D\left( k+1|k\right) 
\text{ ,}
\end{equation*}%
we then have, referring to (\ref{2r_1}) and (\ref{2t}),%
\begin{equation*}
\overset{\cdot }{D}=\overline{F}\,\overline{D}
\end{equation*}%
\begin{equation*}
\overline{D}=\overset{\cdot }{D}+\overline{C}
\end{equation*}%
then%
\begin{equation*}
\overset{\cdot }{D}=\overline{D}-\overline{C}\text{ .}
\end{equation*}%
Continuing,%
\begin{equation*}
\overline{F}\,\overline{D}=\overline{D}-\overline{C}
\end{equation*}%
\begin{equation*}
\overline{F}\,\overline{D}-\overline{D}=\left( \overline{F}-I\right) 
\overline{D}=-\overline{C}\text{ ;}
\end{equation*}%
hence,%
\begin{equation*}
\overline{D}=-\left( \overline{F}-I\right) ^{-1}\overline{C}\text{ .}
\end{equation*}%
Using (\ref{2al}), and (\ref{2am}) 
\begin{equation*}
\overline{D}=-\left( \left[ 
\begin{array}{cc}
1-\alpha & \left( 1-\alpha \right) T \\ 
-\beta /T & 1-\beta%
\end{array}%
\right] -\left[ 
\begin{array}{cc}
1 & 0 \\ 
0 & 1%
\end{array}%
\right] \right) ^{-1}\left( -\left[ 
\begin{array}{c}
\alpha \\ 
\beta /T%
\end{array}%
\right] \right)
\end{equation*}%
\begin{equation*}
=\left[ 
\begin{array}{cc}
-\alpha & \left( 1-\alpha \right) T \\ 
-\beta /T & -\beta%
\end{array}%
\right] ^{-1}\left[ 
\begin{array}{c}
\alpha \\ 
\beta /T%
\end{array}%
\right]
\end{equation*}%
\begin{equation*}
=\frac{\left[ 
\begin{array}{cc}
-\beta & -\left( 1-\alpha \right) T \\ 
\beta /T & -\alpha%
\end{array}%
\right] }{\beta }\left[ 
\begin{array}{c}
\alpha \\ 
\beta /T%
\end{array}%
\right] =\left[ 
\begin{array}{c}
-1 \\ 
0%
\end{array}%
\right] \text{ .}
\end{equation*}%
Hence,%
\begin{equation*}
\overset{\cdot }{D}=\overline{F}\,\overline{D}=\left[ 
\begin{array}{cc}
1-\alpha & \left( 1-\alpha \right) T \\ 
-\beta /T & 1-\beta%
\end{array}%
\right] \left[ 
\begin{array}{c}
-1 \\ 
0%
\end{array}%
\right] =\left[ 
\begin{array}{c}
\alpha -1 \\ 
\beta /T%
\end{array}%
\right] \text{ .}
\end{equation*}%
Finally, the steady-state time updated total covariance is obtained by
substituting into (\ref{2u}), 
\begin{equation*}
\overset{\cdot }{S}=\left[ 
\begin{array}{cc}
\overset{\cdot }{S}_{11} & \overset{\cdot }{S}_{12} \\ 
\overset{\cdot }{S}_{21} & \overset{\cdot }{S}_{22}%
\end{array}%
\right] =\overset{\cdot }{M}+\Phi \overline{D}\Lambda \overline{D}^{\prime
}\Phi ^{\prime }
\end{equation*}%
\begin{equation*}
=\left[ 
\begin{array}{cc}
\overset{\cdot }{M}_{11} & \overset{\cdot }{M}_{12} \\ 
\overset{\cdot }{M}_{21} & \overset{\cdot }{M}_{22}%
\end{array}%
\right] +\left[ 
\begin{array}{cc}
1 & T \\ 
0 & 1%
\end{array}%
\right] \left[ 
\begin{array}{c}
-1 \\ 
0%
\end{array}%
\right] \cdot \Lambda \cdot \left[ 
\begin{array}{cc}
-1 & 0%
\end{array}%
\right] \left[ 
\begin{array}{cc}
1 & 0 \\ 
T & 1%
\end{array}%
\right]
\end{equation*}%
\begin{equation*}
=\left[ 
\begin{array}{cc}
\overset{\cdot }{M}_{11} & \overset{\cdot }{M}_{12} \\ 
\overset{\cdot }{M}_{21} & \overset{\cdot }{M}_{22}%
\end{array}%
\right] +\allowbreak \left[ 
\begin{array}{cc}
\Lambda & 0 \\ 
0 & 0%
\end{array}%
\right] \text{ .}
\end{equation*}%
So,%
\begin{equation*}
\left[ 
\begin{array}{cc}
\overset{\cdot }{S}_{11} & \overset{\cdot }{S}_{12} \\ 
\overset{\cdot }{S}_{21} & \overset{\cdot }{S}_{22}%
\end{array}%
\right] =\frac{N}{\alpha \left( 4-2\alpha -\beta \right) }\left[ 
\begin{array}{cc}
2\alpha ^{2}+2\beta +\alpha \beta & \beta \left( 2\alpha +\beta \right) /T
\\ 
\beta \left( 2\alpha +\beta \right) /T & 2\beta ^{2}/T^{2}%
\end{array}%
\right]
\end{equation*}%
\begin{equation*}
+\frac{q_{22}}{\left( -4\alpha \beta +\alpha \beta ^{2}+2\alpha ^{2}\beta
\right) }\left[ 
\begin{array}{cc}
T^{2}\left( -2+\alpha \right) & T\left( -2\alpha -\beta +\alpha \beta
+\alpha ^{2}\right) \\ 
T\left( -2\alpha -\beta +\alpha \beta +\alpha ^{2}\right) & \left( -2\beta
+2\alpha \beta -2\alpha ^{2}+\alpha ^{3}\right)%
\end{array}%
\right]
\end{equation*}%
\begin{equation*}
+\left[ 
\begin{array}{cc}
0 & 0 \\ 
0 & q_{22}%
\end{array}%
\right] \text{ }+\allowbreak \left[ 
\begin{array}{cc}
\Lambda & 0 \\ 
0 & 0%
\end{array}%
\right] \text{ .}
\end{equation*}

In particular,%
\begin{equation*}
\overset{\cdot }{S}_{11}=\overset{\cdot }{M}_{11}+\Lambda
\end{equation*}%
\begin{equation}
=\frac{N\left( 2\alpha ^{2}+2\beta +\alpha \beta \right) }{\alpha \left(
4-2\alpha -\beta \right) }+\frac{q_{22}T^{2}\left( -2+\alpha \right) }{%
\left( -4\alpha \beta +\alpha \beta ^{2}+2\alpha ^{2}\beta \right) }+\Lambda
\label{2ap}
\end{equation}%
\begin{equation*}
\overset{\cdot }{S}_{21}=\overset{\cdot }{M}_{21}
\end{equation*}%
\begin{equation}
=\frac{N\beta \left( 2\alpha +\beta \right) }{\alpha \left( 4-2\alpha -\beta
\right) T}+\frac{q_{22}T\left( -2\alpha -\beta +\alpha \beta +\alpha
^{2}\right) }{\left( -4\alpha \beta +\alpha \beta ^{2}+2\alpha ^{2}\beta
\right) }\text{ .}  \label{2aq}
\end{equation}%
We turn our attention to (\ref{2x_1}). In steady-state,%
\begin{equation*}
\widetilde{H}S\left( k+1|k\right) \widetilde{H}^{\prime }+\widetilde{N}=%
\left[ 
\begin{array}{cc}
1 & 0%
\end{array}%
\right] \left[ 
\begin{array}{cc}
\overset{\cdot }{S}_{11} & \overset{\cdot }{S}_{12} \\ 
\overset{\cdot }{S}_{21} & \overset{\cdot }{S}_{22}%
\end{array}%
\right] \left[ 
\begin{array}{c}
1 \\ 
0%
\end{array}%
\right] +N+\Lambda
\end{equation*}%
\begin{equation}
=\overset{\cdot }{S}_{11}+N+\Lambda \text{ .}  \label{2ar}
\end{equation}%
Also,%
\begin{equation}
\widetilde{H}\Phi D\left( k|k\right) \Lambda \left( W\frac{\partial u}{%
\partial \lambda }\right) ^{\prime }=\left[ 
\begin{array}{cc}
1 & 0%
\end{array}%
\right] \left[ 
\begin{array}{cc}
1 & T \\ 
0 & 1%
\end{array}%
\right] \left[ 
\begin{array}{c}
-1 \\ 
0%
\end{array}%
\right] \cdot \Lambda \cdot 1\cdot 1=-\Lambda \text{ ,}  \label{2as}
\end{equation}%
and%
\begin{equation}
\left( W\frac{\partial u}{\partial \lambda }\right) \Lambda D\left(
k|k\right) ^{\prime }\Phi ^{\prime }\widetilde{H}^{\prime }=-\Lambda \text{ ,%
}  \label{2at}
\end{equation}%
and%
\begin{equation}
S\left( k+1|k\right) \widetilde{H}^{\prime }=\left[ 
\begin{array}{cc}
\overset{\cdot }{S}_{11} & \overset{\cdot }{S}_{12} \\ 
\overset{\cdot }{S}_{12} & \overset{\cdot }{S}_{22}%
\end{array}%
\right] \left[ 
\begin{array}{c}
1 \\ 
0%
\end{array}%
\right] =\left[ 
\begin{array}{c}
\overset{\cdot }{S}_{11} \\ 
\overset{\cdot }{S}_{12}%
\end{array}%
\right] \text{ ,}  \label{2au}
\end{equation}%
and%
\begin{equation}
\Phi D\left( k|k\right) \Lambda \left( W\frac{\partial u}{\partial \lambda }%
\right) ^{\prime }=\left[ 
\begin{array}{cc}
1 & T \\ 
0 & 1%
\end{array}%
\right] \left[ 
\begin{array}{c}
-1 \\ 
0%
\end{array}%
\right] \cdot \Lambda \cdot 1\cdot 1=\left[ 
\begin{array}{c}
-\Lambda \\ 
0%
\end{array}%
\right] \text{ .}  \label{2av}
\end{equation}

Substituting (\ref{2ar}), (\ref{2as}), (\ref{2at}), (\ref{2au}) and (\ref%
{2av}) into (\ref{2x_1}) gives%
\begin{equation*}
\left[ 
\begin{array}{c}
\alpha \\ 
\beta /T%
\end{array}%
\right] \left( \overset{\cdot }{S}_{11}+N+\Lambda -\Lambda -\Lambda \right) =%
\left[ 
\begin{array}{c}
\overset{\cdot }{S}_{11} \\ 
\overset{\cdot }{S}_{12}%
\end{array}%
\right] +\left[ 
\begin{array}{c}
-\Lambda \\ 
0%
\end{array}%
\right] \text{ .}
\end{equation*}%
This, written as two scalar equations,%
\begin{equation*}
\alpha \left( \overset{\cdot }{S}_{11}+N+\Lambda -\Lambda -\Lambda \right)
=\alpha \left( \overset{\cdot }{S}_{11}+N-\Lambda \right) =\overset{\cdot }{S%
}_{11}-\Lambda
\end{equation*}%
\begin{equation*}
\frac{\beta }{T}\left( \overset{\cdot }{S}_{11}+N+\Lambda -\Lambda -\Lambda
\right) =\frac{\beta }{T}\left( \overset{\cdot }{S}_{11}+N+-\Lambda \right) =%
\overset{\cdot }{S}_{12}\text{ .}
\end{equation*}%
Substituting (\ref{2ap}) and (\ref{2aq})%
\begin{equation*}
\alpha \left( \left( \frac{N\left( 2\alpha ^{2}+2\beta +\alpha \beta \right) 
}{\alpha \left( 4-2\alpha -\beta \right) }+\frac{q_{22}T^{2}\left( -2+\alpha
\right) }{\left( -4\alpha \beta +\alpha \beta ^{2}+2\alpha ^{2}\beta \right) 
}+\Lambda \right) +N-\Lambda \right)
\end{equation*}%
\begin{equation*}
=\frac{N\left( 2\alpha ^{2}+2\beta +\alpha \beta \right) }{\alpha \left(
4-2\alpha -\beta \right) }+\frac{q_{22}T^{2}\left( -2+\alpha \right) }{%
\left( -4\alpha \beta +\alpha \beta ^{2}+2\alpha ^{2}\beta \right) }+\Lambda
-\Lambda
\end{equation*}%
\begin{equation*}
\beta \left( \left( \frac{N\left( 2\alpha ^{2}+2\beta +\alpha \beta \right) 
}{\alpha \left( 4-2\alpha -\beta \right) }+\frac{q_{22}T^{2}\left( -2+\alpha
\right) }{\left( -4\alpha \beta +\alpha \beta ^{2}+2\alpha ^{2}\beta \right) 
}+\Lambda \right) +N+-\Lambda \right)
\end{equation*}%
\begin{equation*}
=\left( \frac{N\beta \left( 2\alpha +\beta \right) }{\alpha \left( 4-2\alpha
-\beta \right) T}+\frac{q_{22}T\left( -2\alpha -\beta +\alpha \beta +\alpha
^{2}\right) }{\left( -4\alpha \beta +\alpha \beta ^{2}+2\alpha ^{2}\beta
\right) }\right) \cdot T
\end{equation*}%
\begin{equation*}
\alpha \left( \frac{N\left( 2\alpha ^{2}+2\beta +\alpha \beta \right) }{%
\alpha \left( 4-2\alpha -\beta \right) }+\frac{q_{22}T^{2}\left( -2+\alpha
\right) }{\left( -4\alpha \beta +\alpha \beta ^{2}+2\alpha ^{2}\beta \right) 
}+N\right) =\frac{N\left( 2\alpha ^{2}+2\beta +\alpha \beta \right) }{\alpha
\left( 4-2\alpha -\beta \right) }+\frac{q_{22}T^{2}\left( -2+\alpha \right) 
}{\left( -4\alpha \beta +\alpha \beta ^{2}+2\alpha ^{2}\beta \right) }
\end{equation*}%
\begin{equation*}
\beta \left( \frac{N\left( 2\alpha ^{2}+2\beta +\alpha \beta \right) }{%
\alpha \left( 4-2\alpha -\beta \right) }+\frac{q_{22}T^{2}\left( -2+\alpha
\right) }{\left( -4\alpha \beta +\alpha \beta ^{2}+2\alpha ^{2}\beta \right) 
}+N\right) =\frac{N\beta \left( 2\alpha +\beta \right) }{\alpha \left(
4-2\alpha -\beta \right) }+\frac{q_{22}T^{2}\left( -2\alpha -\beta +\alpha
\beta +\alpha ^{2}\right) }{\left( -4\alpha \beta +\alpha \beta ^{2}+2\alpha
^{2}\beta \right) }\text{ .}
\end{equation*}%
We define $\rho =q_{22}T^{2}/N$. With $q_{22}$ in $\left( m/\sec \right)
^{2} $, $T$ in seconds and $N$ in $m^{2}$, $\rho $ is unitless. Substituting
this in the previous two equations, we obtain%
\begin{equation*}
\alpha \left( \frac{\left( 2\alpha ^{2}+2\beta +\alpha \beta \right) }{%
\alpha \left( 4-2\alpha -\beta \right) }+\frac{\rho \left( -2+\alpha \right) 
}{\left( -4\alpha \beta +\alpha \beta ^{2}+2\alpha ^{2}\beta \right) }%
+1\right) =\frac{\left( 2\alpha ^{2}+2\beta +\alpha \beta \right) }{\alpha
\left( 4-2\alpha -\beta \right) }+\frac{\rho \left( -2+\alpha \right) }{%
\left( -4\alpha \beta +\alpha \beta ^{2}+2\alpha ^{2}\beta \right) }
\end{equation*}%
\begin{equation*}
\beta \left( \frac{\left( 2\alpha ^{2}+2\beta +\alpha \beta \right) }{\alpha
\left( 4-2\alpha -\beta \right) }+\frac{\rho \left( -2+\alpha \right) }{%
\left( -4\alpha \beta +\alpha \beta ^{2}+2\alpha ^{2}\beta \right) }%
+1\right) =\frac{\beta \left( 2\alpha +\beta \right) }{\alpha \left(
4-2\alpha -\beta \right) }+\frac{\rho \left( -2\alpha -\beta +\alpha \beta
+\alpha ^{2}\right) }{\left( -4\alpha \beta +\alpha \beta ^{2}+2\alpha
^{2}\beta \right) }\text{ .}
\end{equation*}%
We divide the previous two equations and cancel an $\alpha $,%
\begin{equation*}
\frac{\alpha }{\beta }=\frac{\left( 2\alpha ^{2}+2\beta +\alpha \beta
\right) \left( -4\beta +\beta ^{2}+2\alpha \beta \right) +\rho \left(
-2+\alpha \right) \left( 4-2\alpha -\beta \right) }{\beta \left( 2\alpha
+\beta \right) \left( -4\beta +\beta ^{2}+2\alpha \beta \right) +\rho \left(
-2\alpha -\beta +\alpha \beta +\alpha ^{2}\right) \left( 4-2\alpha -\beta
\right) }\text{ .}
\end{equation*}%
$\bigskip $Cross multiplying gives:%
\begin{equation}
\allowbreak \allowbreak 2\beta ^{4}+\left( 4\alpha -8\right) \beta ^{3}+\rho
\left( \allowbreak \left( \alpha ^{2}-2\alpha +2\right) \beta ^{2}+\left(
3\alpha ^{3}-10\alpha ^{2}\allowbreak +12\alpha -8\right) \beta +\left(
2\alpha ^{4}-8\alpha ^{3}+8\alpha ^{2}\right) \right) =0\text{ .}
\label{2aw}
\end{equation}%
Equation (\ref{2aw}) gives our relationship between $\alpha $ and $\beta $.
The noise ratio $\rho $ is a parameter in the equation which is known, or at
least known to be within a range. Equation (\ref{2aw}) may be factored%
\begin{equation}
\left( \beta +\left( 2\alpha -4\right) \right) \left( 2\beta ^{3}+\rho
\left( \left( \alpha ^{2}-2\alpha +2\right) \beta +\alpha ^{2}\left( \alpha
-2\right) \right) \right) =0\text{ .}  \label{2ax}
\end{equation}

Hence, $\beta =(4-2\alpha )$ is a solution that is independent of $\rho $.
This solution is not permitted however since it violates condition 3. of (%
\ref{2am4}). There is a second real solution for $\beta $ given $\alpha $
(which depends on $\rho $.) The remaining two solutions for $\beta $ given $%
a $ may be a complex conjugate pair. The table below gives some
representative solutions to (\ref{2aw}). The two real solutions are
presented. The second one listed we don't use because of the condition 3.
The Newton-Raphson method may be used to compute all of the solutions to (%
\ref{2aw}).

\begin{center}
\begin{tabular}{|ccc|}
\hline
$\rho $ & $\alpha $ & $\beta $ \\ \hline
2 & 0.2 & \multicolumn{1}{l|}{0.04385, 3.6} \\ 
4 & 0.2 & \multicolumn{1}{l|}{0.04386, 3.6} \\ 
6 & 0.2 & \multicolumn{1}{l|}{0.04389, 3.6} \\ 
6 & 0.4 & \multicolumn{1}{l|}{0.1866, 3.2} \\ 
8 & 0.2 & \multicolumn{1}{l|}{0.04389, 3.6} \\ 
8 & 0.4 & \multicolumn{1}{l|}{0.1870, 3.2} \\ 
10 & 0.2 & \multicolumn{1}{l|}{0.04389, 3.6} \\ 
10 & 0.4 & \multicolumn{1}{l|}{0.1873, 3.2} \\ 
10 & 0.5 & \multicolumn{1}{l|}{0.2959, 3.0} \\ \hline
\end{tabular}
\end{center}

These $\alpha $ and $\beta $ give that the eigenvalues of $\overline{F}$, in
(\ref{2al1}), have norm less than $1$. Hence, by Theorem 2.1, page 64 of 
\cite{AandM}, $M_{N}$ and $M_{Q}$, the solutions to (\ref{2am1i}) and (\ref%
{2am1ii}) respectively, exist are unique and are positive definite.

\section{Summary and Conclusions}

In this paper we considered some topics in radar sensor bias. We presented
an algorithm that estimates the absolute bias of two sensors when the
relative bias between the sensors is given. The algorithm uses the relative
bias, which is given in rectangular coordinates, as a constraint. The
absolute biases, in spherical coordinates, for the sensors are obtained by
the solution to an optimization problem that exploits the
spherical-to-rectangular coordinate conversion. We presented a reduced-state
filter that is designed for performance with sensor bias. The filter is
reduced-state since it does not contain additional bias states. The filter
design is influenced by the filter in \cite{MooRief}. It may be viewed as a
dual design (in the control theory sense) to the \ filter in \cite{MooRief}.

A flow diagram for processing radar data with bias may contain these stages:

1. Estimate state with the $\alpha -\beta $ filter optimized for measurement
bias, as presented in Section 3.

2. For a multi-sensor problem, estimate the relative sensor bias using an
optimized algorithm such as in \cite{Mark}.

3. Continue by estimating the absolute bias for each sensor using the
algorithm presented in Section 2.

\section{Appendix: Transformation from $ENU(1)$ to $ENU(2)$}

In this appendix we present the transformation from the $ENU(1)$ to the $%
ENU(2)$ coordinate systems. But first, consider the transformation from $ECI$
to $ENU$. Consider an $ENU$ coordinate axis located at longitude-latitude $%
\Omega -L$ and define the rotation matrix%
\begin{equation*}
T_{ECI2ENU}=\left[ 
\begin{array}{ccc}
1 & 0 & 0 \\ 
0 & \cos L & -\sin L \\ 
0 & \sin L & \cos L%
\end{array}%
\right] \left[ 
\begin{array}{ccc}
\cos \Omega & 0 & -\sin \Omega \\ 
0 & 1 & 0 \\ 
\sin \Omega & 0 & \cos \Omega%
\end{array}%
\right] \left[ 
\begin{array}{ccc}
0 & 1 & 0 \\ 
0 & 0 & 1 \\ 
1 & 0 & 0%
\end{array}%
\right]
\end{equation*}

%
\begin{equation*}
=\left[ 
\begin{array}{ccc}
-\sin \Omega & \cos \Omega & 0 \\ 
-\sin L\cos \Omega & -\sin L\sin \Omega & \cos L \\ 
\cos L\cos \Omega & \cos L\sin \Omega & \sin L%
\end{array}%
\right] \allowbreak
\end{equation*}%
with $T_{ENU2ECI}=T_{ECI2ENU}^{\prime }$. For position, we need to include a
translation, so that for a given position vector in $ENU$ coordinates%
\begin{equation*}
P_{ECI}=\frac{r_{ee}}{\sqrt{1-e^{2}\sin ^{2}L}}\left[ 
\begin{array}{c}
\cos L\cos \Omega \\ 
\cos L\sin \Omega \\ 
\left( 1-e^{2}\right) \sin L%
\end{array}%
\right] +T_{ENU2ECI}P_{ENU}
\end{equation*}%
where $r_{ee}$ is the earth's equatorial radius and $e$ is the earth's
eccentricity. For velocity, use the rotation alone.

Let the $ENU(1)$, $ENU(2)$ coordinate system be located at
longitude-latitude $\Omega _{1}-L_{1}$ and $\Omega _{2}-L_{2}$ respectively.
Next we consider our transformation going from $\Omega _{1}-L_{1}$ to $%
\Omega _{2}-L_{2}$. The rotation part of this transformation can be
represented by the matrix below (going from $ENU(1)$ to $ENU(2)$). There are 
$3$ steps. First, the $ENU(1)$ coordinates are rotated down to the equator.
Second, these coordinates are rotated along the equator by the longitude
difference. Third is the rotation up to the latitude of the $ENU(2)$ system.%
\begin{equation*}
T_{ENU(1)2ENU(2)}=\left[ 
\begin{array}{ccc}
1 & 0 & 0 \\ 
0 & \cos L_{2} & -\sin L_{2} \\ 
0 & \sin L_{2} & \cos L_{2}%
\end{array}%
\right] \left[ 
\begin{array}{ccc}
\cos \left( \Omega _{2}-\Omega _{1}\right) & 0 & -\sin \left( \Omega
_{2}-\Omega _{1}\right) \\ 
0 & 1 & 0 \\ 
\sin \left( \Omega _{2}-\Omega _{1}\right) & 0 & \cos \left( \Omega
_{2}-\Omega _{1}\right)%
\end{array}%
\right] \left[ 
\begin{array}{ccc}
1 & 0 & 0 \\ 
0 & \cos L_{1} & \sin L_{1} \\ 
0 & -\sin L_{1} & \cos L_{1}%
\end{array}%
\right]
\end{equation*}%
\begin{equation*}
=\left[ 
\begin{array}{ccc}
\cos \left( \Omega _{2}-\Omega _{1}\right) & \sin L_{1}\sin \left( \Omega
_{2}-\Omega _{1}\right) & -\cos L_{1}\sin \left( \Omega _{2}-\Omega
_{1}\right) \\ 
-\sin L_{2}\sin \left( \Omega _{2}-\Omega _{1}\right) & \cos L_{1}\cos
L_{2}+\sin L_{1}\sin L_{2}\cos \left( \Omega _{2}-\Omega _{1}\right) & \cos
L_{2}\sin L_{1}-\cos L_{1}\sin L_{2}\cos \left( \Omega _{2}-\Omega
_{1}\right) \\ 
\cos L_{2}\sin \left( \Omega _{2}-\Omega _{1}\right) & \cos L_{1}\sin
L_{2}-\cos L_{2}\sin L_{1}\cos \left( \Omega _{2}-\Omega _{1}\right) & \sin
L_{1}\sin L_{2}+\cos L_{1}\cos L_{2}\cos \left( \Omega _{2}-\Omega
_{1}\right)%
\end{array}%
\right] \allowbreak
\end{equation*}%
(Note%
\begin{equation*}
\left[ 
\begin{array}{ccc}
\cos \left( \Omega _{2}-\Omega _{1}\right) & 0 & -\sin \left( \Omega
_{2}-\Omega _{1}\right) \\ 
0 & 1 & 0 \\ 
\sin \left( \Omega _{2}-\Omega _{1}\right) & 0 & \cos \left( \Omega
_{2}-\Omega _{1}\right)%
\end{array}%
\right] =\left[ 
\begin{array}{ccc}
\cos \Omega _{1}\cos \Omega _{2}+\sin \Omega _{1}\sin \Omega _{2} & 0 & 
-\cos \Omega _{1}\sin \Omega _{2}+\cos \Omega _{2}\sin \Omega _{1} \\ 
0 & 1 & 0 \\ 
\cos \Omega _{1}\sin \Omega _{2}-\cos \Omega _{2}\sin \Omega _{1} & 0 & \cos
\Omega _{1}\cos \Omega _{2}+\sin \Omega _{1}\sin \Omega _{2}%
\end{array}%
\right] \allowbreak \allowbreak
\end{equation*}%
\begin{equation*}
=\left[ 
\begin{array}{ccc}
\cos \Omega _{2} & 0 & -\sin \Omega _{2} \\ 
0 & 1 & 0 \\ 
\sin \Omega _{2} & 0 & \cos \Omega _{2}%
\end{array}%
\right] \left[ 
\begin{array}{ccc}
\cos \Omega _{1} & 0 & \sin \Omega _{1} \\ 
0 & 1 & 0 \\ 
-\sin \Omega _{1} & 0 & \cos \Omega _{1}%
\end{array}%
\right]
\end{equation*}%
so that the rotation $T_{ENU(1)2ENU(2)}$ is a rotation from the first
coordinates down to $ECI$ and then up to the second coordinates.) We have%
\begin{equation*}
T_{ENU(2)2ENU(1)}=T_{ENU(1)2ENU(2)}^{\prime }
\end{equation*}

The position vector from the $ENU(1)$ to the $ENU(2)$ coordinate axes (in $%
ECI$ coordinates) is%
\begin{equation*}
P_{ENU(1)2ENU(2),ECI}=r_{ee}\frac{\left[ 
\begin{array}{c}
\cos L_{2}\cos \Omega _{2} \\ 
\cos L_{2}\sin \Omega _{2} \\ 
\left( 1-e^{2}\right) \sin L_{2}%
\end{array}%
\right] }{\sqrt{1-e^{2}\sin ^{2}L_{2}}}-r_{ee}\frac{\left[ 
\begin{array}{c}
\cos L_{1}\cos \Omega _{1} \\ 
\cos L_{1}\sin \Omega _{1} \\ 
\left( 1-e^{2}\right) \sin L_{1}%
\end{array}%
\right] }{\sqrt{1-e^{2}\sin ^{2}L_{1}}}
\end{equation*}%
and in the other coordinates this vector is%
\begin{equation*}
P_{ENU(1)2ENU(2),ENU\left( i\right) }=T_{ECI2ENU(i)}P_{ENU(1)2ENU(2),ECI}
\end{equation*}%
The total position coordinate transformation, including translation can be
represented by%
\begin{equation*}
P_{ENU(2)}=-P_{ENU(1)2ENU(2),ENU\left( 2\right) }+T_{ENU(1)2ENU(2)}P_{ENU(1)}
\end{equation*}%
The total velocity coordinate transformation is given by the rotation alone.

\end{document}